\documentclass[11pt]{amsart}

\usepackage{amssymb}
\usepackage{microtype}
\numberwithin{equation}{section}
\usepackage{xcolor}
\usepackage{hyperref}
\hypersetup{   
colorlinks=true,       
linkcolor=black,          
citecolor=black,        
filecolor=black,      
urlcolor=black}

\newtheorem{theorem}{Theorem}[section]

\newtheorem{lemma}[theorem]{Lemma}
\newtheorem{lem}[theorem]{Lemma}
\newtheorem{prop}[theorem]{Proposition}
\newtheorem{proposition}[theorem]{Proposition}

\newtheorem{corollary}[theorem]{Corollary}

\newcommand{\N}{\mathbb{N}}
\newcommand{\R}{\mathbb{R}}
\newcommand{\C}{\mathbb{C}}
\newcommand{\PSH}{{\rm PSH}}

\theoremstyle{definition}
\newtheorem{definition}[theorem]{Definition}

\newtheorem{example}[theorem]{Example}

\newtheorem{remark}[theorem]{Remark}
\newtheorem*{ackn}{Acknowledgements}

 \subjclass[2010]{31C45, 2U15, 32U40,
32W20, 35J66, 35J96}

\keywords{Complex Monge-Amp\`ere operator, Dirichlet eigenvalue, iterative
method, uniqueness, plurisubharmonic envelope}

\begin{document}

\title[A new approach to the Monge-Amp\`ere eigenvalue problem]{A new approach to the Monge-Amp\`ere eigenvalue problem}

\author{Chinh H. Lu}
\address{Chinh H. Lu: Institut Universitaire de France \& Univ Angers, CNRS, LAREMA, SFR MATHSTIC, F-49000 Angers, France}
\email{\href{hoangchinh.lu@univ-angers.fr}{hoangchinh.lu@univ-angers.fr}} 

\author{Ahmed Zeriahi}
\address{Ahmed Zeriahi: Institut de Math\'ematiques de Toulouse   \\ Universit\'e de Toulouse \\
118 route de Narbonne \\
31400 Toulouse, France\\}

\email{\href{mailto:ahmed.zeriahi@math.univ-toulouse.fr}{ahmed.zeriahi@math.univ-toulouse.fr}}

 \begin{abstract}
 	We study the eigenvalue problem for the complex Monge-Amp\`ere operator in bounded hyperconvex domains in $\mathbb C^n$, where the right-hand side is a non-pluripolar positive Borel measure. We establish the uniqueness of eigenfunctions in the finite energy class introduced by Cegrell, up to positive multiplicative constants, and provide a Rayleigh quotient type formula for computing the eigenvalue. 

Under a natural continuity assumption on the measure, we further show that both the eigenvalue and eigenfunctions can be obtained via an iterative procedure starting from any negative finite energy function. 

Our approach relies on the fine properties of plurisubharmonic envelopes, which allow a partial sublinearization of the nonlinear problem. As far as we know, this method is new, even in the linear case, and not only yields new results but also significantly simplifies existing arguments in the literature. Moreover, it extends naturally to the setting of complex Hessian operators. 

Finally, by translating our results from the complex Monge-Amp\`ere setting via a logarithmic transformation, we also obtain several interesting analogues for the real Monge-Amp\`ere operator. 
\end{abstract}

\date{\today}
 
 \maketitle
 
{\small \tableofcontents}

\section{Introduction}

It is well known that a linear elliptic operator admits infinitely many positive eigenvalues. Among these, the first eigenvalue is of particular interest: it is the only one associated with a negative eigenfunction, it admits a variational characterization involving the Dirichlet energy, and its corresponding eigenspace is one-dimensional.

In a foundational paper, P.-L. Lions \cite{Lio85} investigated the eigenvalue problem for the real Monge-Amp\`ere operator. For smooth data, he established the existence of a unique eigenvalue and showed that all convex eigenfunctions are proportional, thus providing a satisfactory extension of the classical one-dimensional theory. Lions further conjectured that his approach could be extended to the complex Monge-Amp\`ere operator. This conjecture was recently confirmed by Badiane and the second author in \cite{BZ23}. The a priori estimates in \cite{BZ23} are not straightforward consequences of Lions' original arguments. While the aforementioned works assume smooth domains and densities, it is also possible to consider degenerate settings, as demonstrated in \cite{BZ24} and \cite{LeQN18}.

  Let $\Omega$ be a bounded hyperconvex domain in $\mathbb C^n$, and let $\mu$ be a positive Borel measure on $\Omega$ (with locally finite mass). The eigenvalue problem consists of finding a pair $(\lambda,u)$, where $\lambda >0$ and $u\not \equiv  0$ lies in the Cegrell finite energy class  $\mathcal E^1(\Omega)$, whose definition will be recalled in the next section, such that the following equation holds
\begin{equation} \label{eq: MA eigen}
	(dd^c u)^n = (-\lambda u)^n \mu
\end{equation} 
in the weak sense of Radon measures. If $(\lambda,u)$ solves \eqref{eq: MA eigen} then $\lambda$ is called an eigenvalue while $u$ is called an eigenfunction. 

As shown in \cite[Theorem 5.2]{CLM24}, building on an argument of L\^e \cite{LeQN18}, the eigenvalue is uniquely determined if it exists. In \cite{CLM24}, the authors  assumed the data are smooth but the argument also works for the general case.  However, the uniqueness of eigenfunctions up to positive multiplicative constants is significantly more delicate. It was established in \cite{BZ23}, and adapted to the Hessian case in \cite{CLM24}, under strong regularity assumptions on the data: the boundary of $\Omega$ is smooth and $\mu$ has positive smooth density. The key idea there is to linearize the equation  using a smooth solution, which in turn requires the existence of the latter, hence the need for smooth data.

Our first main result establishes the uniqueness of eigenfunctions, and provides a  Rayleigh quotient type formula  computing  the eigenvalue, merely assuming  that $\mu$ vanishes on pluripolar sets. 

\begin{theorem} \label{thm: main unique}
 Assume $\mu$ is a non-pluripolar positive Borel measure in $\Omega$, and $(\lambda,\varphi)$ solves the eigenvalue problem \eqref{eq: MA eigen}. Then  $\lambda=\lambda_1(\mu)$ is given by 
 	\[
 	\lambda_1(\mu) = \inf \left \{ \frac{E(u)}{I_{\mu}(u)} \; , \; u \in \mathcal E^1(\Omega) \setminus \{0\} \right \}. 
 	\]
  Moreover, if $\psi \in \mathcal E^1(\Omega) \setminus \{0\}$ satisfies $(dd^c \psi)^n \geq (-\lambda \psi)^n \mu$, then $(\lambda,\psi)$ actually solves \eqref{eq: MA eigen}, and $\psi=c \varphi$, for some positive constant $c$. 
 \end{theorem}

The space $\mathcal E^1(\Omega)$  appearing in the above theorem is the Cegrell space of plurisubharmonic functions in $\Omega$ with finite Monge-Ampère energy (see Section 2) and the 
functionals appearing are defined on the space $\mathcal E^1(\Omega)$ as follows
\[
I_{\mu}:=\int_{\Omega} (-u)^{n+1} d\mu\; , \; E(u):=\int_{\Omega} (-u)(dd^c u)^n.
\]  
In the real setting, the variational formulation computing $\lambda_1(\mu)$ was proved by Tso \cite{Tso90}.  The main novelty of our argument is that instead of using a smooth solution to linearize the equation, we use plurisubharmonic envelopes, as in \cite{Sal25}, \cite{ACLR24,ACLR25}.

Let us stress that in the above uniqueness result we do not assume any regularity assumption on the non-pluripolar positive measure $\mu$. To solve \eqref{eq: MA eigen}, however, some regularity condition on $\mu$ will be needed, as we now describe. 
We can use the strategy introduced by Lions in \cite{Lio85} which  proceeds as follows. We let $\gamma_1$ be the supremum of all positive $\gamma$ for which the equation 
\[
(dd^c u)^n = (1 -\gamma u)^n \mu 
\]
has a solution. We then consider a sequence of solutions $(\gamma_j, u_j)$ with $\gamma_j \nearrow \gamma_1$. As $\|u_j\|_{\infty} \to +\infty$, the rescaled functions $v_j := \|u_j\|_{\infty}^{-1} u_j$ converge (up to a subsequence) to a bounded solution of the eigenvalue problem \eqref{eq: MA eigen}. This approach was used in \cite{BZ23} for smooth data, but we will show later that it applies to any measure $\mu = f dV$ with $f \in L^p$, $p>1$; see Theorem \ref{thm: eigen Hessian bounded}. 

For less regular measure, e.g. $\mu$ does not have a density with respect to $dV$,  we can solve \eqref{eq: MA eigen} using a variational method. This approach was used in \cite{BZ23, BZ24} and it works when the functional $I_{\mu}$ is continuous on each $\mathcal E^1_C(\Omega)$, consisting of functions $u\in \mathcal E^1(\Omega)$ such that $E(u)\leq C$.

Recently, F. Abedin and J. Kitagawa proposed in \cite{AK20} an inverse iteration method  for the real Monge-Amp\`ere eigenvalue problem. They showed that the iterative scheme converges to an eigenfunction if the latter exists. Shortly after, Q.N. L\^e \cite{LeQN20} showed that the iterative scheme of F. Abedin and J. Kitagawa can be started from any non-trivial data provided that an eigenfunction exists.  Very recently, the second author proved in \cite{Zer25} that this method works in the setting of complex Monge-Amp\`ere equations under the assumption on the existence of the eigenvalue and an eigenfunction. Given $u_0\in \mathcal E^1(\Omega)\setminus \{0\}$, we define inductively  $u_{k+1}\in \mathcal E^1(\Omega)$ as 
\[
(dd^c u_{k+1})^n = R(u_k) (-u_k)^n \mu, \; R(u_k):= \frac{E(u_k)}{I_{\mu}(u_k)}. 
\] 
One of the delicate points here lies in proving that the sequence $(u_k)_{k \in \N}$ does not converge to $0$. This can be done by using the existence of a solution, as in \cite{Zer25}, \cite{AK20}, \cite{LeQN20}. 
In our next result, we observe that along this iterative scheme the energy $E(u_k)$ is increasing, thus starting from any negative function $u_0$,  the sequence $(u_k)$ never converges to $0$. 
\begin{theorem}
	\label{thm: iterative intro}
	Assume $I_{\mu}$ is finite on $\mathcal E^1(\Omega)$, and consider the iterative sequence $(u_k)$ defined as above. Then the sequences $(E(u_k))_{k \in \N}$ and $(I_{\mu}(u_k))_{k \in \N}$ are increasing while $(R(u_k))_{k \in \N}$ is decreasing.  Moreover, if $I_{\mu}$ is continuous on each $\mathcal E^1_C(\Omega)$,  then the sequence $(R(u_k),u_k)_{k \in \N}$ converges to $(\lambda_1, u)$ which solves \eqref{eq: MA eigen}.   
\end{theorem}

The strength of the above result is that it produces an efficient method to prove the existence of the eigenvalue and  a corresponding eigenfunction and approximate them. Compared to \cite{LeQN20}, \cite{AK20}, \cite{Zer25},  we do not assume a priori the existence of  the eigenfunctions, our result thus provides a {\it  new method} to solve \eqref{eq: MA eigen}.

As an application of the uniqueness result in Theorem \ref{thm: main unique}, we obtain the following existence and uniqueness of solution for a class of complex Monge-Amp\`ere equations whose right-hand side is not necessarily increasing on the unknown. We consider the following complex Monge-Amp\`ere equation
\begin{equation} \label{eq:DP intro}
	(dd^c \varphi)^n  -  F(\cdot,\varphi)^n \mu = 0, \; \varphi \in \mathcal E^1(\Omega), 
\end{equation}
where $F :  \Omega \times (-\infty,0] \longrightarrow \R^+$ is a Borel function such that 
\begin{itemize}
\item $ f:= F(\cdot,0)  \in L^{n+1}(\Omega,\mu)$, 
\item for $ \mu$-a.e. $z \in \Omega$,  $t \longmapsto F(z,t)$ is continuous in $]-\infty,0]$. 
\end{itemize}

\begin{theorem}  \label{thm:Existence-Uniqueness intro} We assume $\mu(\Omega)<+\infty$, $I_{\mu}$ is finite on $\mathcal E^1(\Omega)$, and there exists a real number $ \lambda_0 $ such that $0 \leq \lambda_0 < \lambda_1$ and 
$$
\frac{\partial F}{\partial t} \geq - \lambda_0 , \, \, \mathrm{in} \, \,  ]-\infty, 0] \times \Omega.
$$ 
Then the Dirichlet problem \eqref{eq:DP} admits a unique solution $\varphi \in \mathcal E^1(\Omega)$. In particular, if $F(z,0)\equiv 0$, then the unique solution is identically $0$.
\end{theorem}

A direct consequence of this theorem is that the eigenvalue $\lambda_1(\Omega,\mu)$ can be computed by Lions' approach. It is the largest $\lambda>0$ for which the equation $(dd^c u)^n = (1-\lambda u)^n \mu$ admits a solution in $\mathcal E^1(\Omega)$. In Theorem \ref{thm:Existence-Uniqueness2}, we actually solve a more general  Dirichlet problem with continuous boundary values. It is worth noting that an existence result for the Dirichlet problem \eqref{eq:DP intro} was obtained in \cite{CK06} under a strong assumption on $F$ which is not satisfied in our case (e.g. it cannot be applied to \eqref{eq:modifiedMA}).

 Our approach works equally for the complex Hessian operator and gives analogous results, as described in Section \ref{sect: Hessian}.  We refer the reader to the recent work \cite{CLM24}, where the eigenvalue problem is studied for smooth data.  
 
 The continuity assumption on $I_{\mu}$ in Theorem \ref{thm: iterative intro}  is also sufficient to run the variational method, as shown in \cite{BZ24}. This condition is satisfied by a large class of  measures including the ones with $L^p$-density, $p>1$.  It is also satisfied if $\mu$ is dominated by the Monge-Amp\`ere measure of a continuous plurisubharmonic function, as shown in Proposition \ref{prop: I cont}.  The latter criterion turns out to be very useful in the real setting since convex functions are automatically continuous. We use it in the last section to obtain analogous results for the real Monge-Amp\`ere operator by relating the real problem to the complex one via a logarithmic transformation. 
 
 Let $D$ be a bounded convex domain in $\mathbb R^n$ and $\nu$ be a positive Borel measure on $D$ with positive mass $\nu(D)>0$. We look for pairs $(\lambda, u)$ solving 
 \begin{equation}
 	\label{eq: eigenvalue convex}
 	 \text{M}_\R(u) = (-\lambda u)^n \nu, \, \, \, u_{\mid \partial D} = 0. 
 \end{equation}
   where $\lambda > 0$ and $u \in {\rm CV} (D) \cap C^0(\bar D)$ is normalized by $\Vert u\Vert_{C^0(\bar D)}  = 1$.  Here, $\text{M}_{\mathbb{R}}$ denotes the real Monge-Amp\`ere operator. 
 \begin{theorem}\label{thm: MA convex}{\; } 
 	\begin{enumerate}
 		\item  The eigenvalue problem above has at most one solution. 
 		\item If $\nu = \text{M}_{\mathbb R}(v)$, for some $v\in {\rm CV}(D) \cap C^0(\bar D)$, then the eigenvalue problem has a unique solution.    
 	\end{enumerate}
 \end{theorem}
In the uniqueness result (1), we do not impose any condition on the positive Borel measure $\nu$. The condition in (2) is satisfied by a broad class of Borel measures, including those that integrate certain negative convex functions, as shown in our next result. 
\begin{theorem}\label{thm: MA convex N intro}
	Assume $\nu$ is a positive Borel measure on $D$ with positive mass $\nu(D)>0$. If there exists $v \in {\rm CV}(D)$ such that $v<0$ and $\int_D (-v) d\nu <+\infty$, then the Monge-Amp\`ere equation ${\rm M}_{\mathbb R}(u) =\nu$ has a unique solution $u\in {\rm CV}_0(D)$.  
\end{theorem}
This theorem can be viewed as a convex analogue of \cite[Proposition 5.2]{Ceg08}, and it significantly generalizes \cite[Proposition 2.1]{Bl03} and \cite[Theorem 1.1]{Har06}, where the assumption $\nu(D)<+\infty$ is required. In particular, the eigenvalue problem \eqref{eq: eigenvalue convex} is solvable for any positive Borel measure with positive finite total mass. On the way proving Theorem \ref{thm: MA convex N intro}, we also show that any convex function in $D$ whose logarithmic transformation is in the Cegrell class $\mathcal N$ is continuous up to the boundary where it takes the value $0$.

Nevertheless, as illustrated in Example \ref{example: mass},  for any continuous convex functions $0\not \equiv u$ vanishing on the boundary, the function $u_{\alpha}:=-(-u)^{\alpha}$, for $\alpha\in (0,1)$, is continuous and convex, yet its real Monge-Amp\`ere measure does not have finite total mass.

 The proof of Theorem \ref{thm: MA convex} follows, to a large extent, from its complex analogue, as the real Monge-Amp\`ere equation can be translated into a complex one via the logarithmic map
\[
(z_1,...,z_n) \mapsto L(z) = (\log|z_1|,...,\log|z_n|).
\]
We emphasize that no smoothness assumptions are made on either $\partial D$ or the measure $\nu$. Our result thus significantly generalizes earlier works by Lions \cite{Lio85} and L\^e \cite{LeQN18}.

Using the iterative approach, it is also possible to replace the condition in (the second statement of Theorem \ref{thm: MA convex} by  $I_{\mu}$ being continuous on each $\mathcal E^1_C(D)$. Here  the set $\mathcal E^1_C(D)$  consists of functions $w\in {\rm CV}(D)\cap C^0(\bar D)$, vanishing on $\partial D$, and $\int_{D} (-w) {\rm M}_{\mathbb R}(w)\leq C$. 

The paper is organized as follows. In Section \ref{sect: preliminaries} we recall the definitions of functions in Cegrell's classes along with fundamental tools in Pluripotential Theory. Section \ref{sect: eigenvalue problem} is devoted to the systematic study of uniqueness of the eigenfunctions, where we prove Theorem \ref{thm: main unique} (see Corollary \ref{cor: unique}, Theorem \ref{thm: subsol is sol} and Theorem \ref{thm: uniqueness of the eigenvalue}), and we also give a sufficient criterion for the existence of eigenfunctions.

 In Section \ref{sect: Dirichlet problem} we study a general Dirichlet problem  with right-hand side non-increasing in the unknown, and we prove Theorem \ref{thm:Existence-Uniqueness intro}. 
 
 In Section \ref{sect: iterative}, we show that the iterative method introduced in \cite{AK20} can be used to solve the eigenvalue problem, proving Theorem \ref{thm: iterative intro}. 
 
 In the last two sections, we extend our findings to the complex Hessian and real Monge-Amp\`ere settings. 

\begin{ackn}
This work is partially supported by Institut Universitaire de France
and the KRIS project of fondation Charles Defforey. Part of this work was done during the visit of the authors at El Jadida University in April 2025. We thank Omar Alehyane for the invitation and the hospitality. 
\end{ackn}

\section{Preliminaries}\label{sect: preliminaries}

\subsection{The Monge-Amp\`ere  energy functional}
Let $\Omega \Subset \C^n$ be a bounded hyperconvex domain i.e. it admits a continuous negative plurisubharmonic exhausion. Let us define an important convex class of singular plurisubharmonic functions in $\Omega$ suitable for the variational approach.  

First we denote by $\PSH(\Omega) \subset L^1_{\rm loc} (\Omega)$ the convex positive cone of pluribubharmonic functions in $\Omega$. Observe that given  $u \in \PSH(\Omega) \cap C^2(\Omega)$, we can consider the associated  continuous $(n,n)$-form on $\Omega$ defined by
$$
(dd^c u)^n = \text{det} \ \left(\frac{\partial^2 u}{\partial z_j \partial \bar  z_k}\right) \beta^n,
$$ 
This form can be viewed as an $(n,n)$-current on $\Omega$ and can be identified to a Borel measure (with density) on $\Omega$, called the (complex) Monge-Amp\`ere  measure of $u$. Here  $\beta := dd^c \vert z\vert^2$ is the standard K\"ahler form on $\C^n$.

When $u \in \PSH(\Omega) \cap L^{\infty}_{ \rm loc} (\Omega)$, a deep result of E. Bedford and B.A. Taylor \cite {BT76} shows that the   Monge-Amp\`ere  measure  $(dd^c u)^n $ of $u$ on $\Omega$ can be defined as the weak limit of the sequence of Monge-Amp\`ere  measures of any decreasing sequence of smooth plurisubhramonic functions converging to $u$ in $\Omega$. 

 Following Urban Cegrell \cite{Ceg98}, we define the class $\mathcal{E}^0 (\Omega)$ as the set of bounded plurisubharmonic functions $\phi$ on $\Omega$ with boundary values $0$ such that $\int_\Omega (dd^c \phi)^n < + \infty$. Then we define 
$\mathcal{E}^1 (\Omega)$ as the set of plurisubharmonic functions $u$ in $\Omega$ such that there exists a decreasing sequence $(u_j)_{j \in \N}$ in the class $\mathcal{E}^0 (\Omega)$  satisfying $u = \lim_j u_j $ in $\Omega$ and $\sup_j \int_\Omega (-u_j) (dd^c u_j)^n < + \infty$. 
It is easy to see from the definition that  $\mathcal{E}^1 (\Omega)$ is a convex positive cone in $L^1_{\rm loc}(\Omega)$.  

A negative psh function $u$ is in the class $\mathcal F(\Omega)$ if there exists a decreasing sequence $(u_j)_{j \in \N}$ in the class $\mathcal{E}^0 (\Omega)$  satisfying $u = \lim_j u_j $ in $\Omega$ and $\sup_j \int_\Omega  (dd^c u_j)^n < + \infty$. The class $\mathcal E(\Omega)$ consists of all negative psh functions $u$ such that for each relatively compact subset $U\Subset \Omega$, there exists $v\in \mathcal F(\Omega)$ with $u=v$ on $U$. 

Let $\Omega_j$ be an increasing sequence of relatively compact open subsets of $\Omega$ such that $\cup \Omega_j=\Omega$. Given $u\in \mathcal E(\Omega)$, we define $u^j$ to be the largest psh function in $\Omega$ that lies below $u$ on $\Omega\setminus \Omega_j$. Then $u^j$ increases almost everywhere to a psh function $\tilde{u}$. We say that $u\in \mathcal N(\Omega)$ if $\tilde{u}=0$.

Cegrell proved that the complex Monge-Amp\`ere operator extends to the class $\mathcal{E}^1 (\Omega)$ and is continuous under monotone limits in $\mathcal{E}^1 (\Omega)$. Moreover if $u \in \mathcal{E}^1 (\Omega)$, then $\int_\Omega (-u) (dd^c u)^n < + \infty$ (see  \cite{Ceg98}). 

The Monge-Amp\`ere energy functional is defined on the space $\mathcal{E}^1 (\Omega)$  by 
\begin{equation} \label{eq:energyF}
E (\phi) :=   \int_\Omega (-\phi) (dd^c \phi)^n, \; \phi \in \mathcal{E}^1 (\Omega).
\end{equation}
 
\begin{lem}   \label{Primitive}  Let $\phi_0,\phi_1 \in \mathcal E^1(\Omega)$ and $\phi_t=(1-t)\phi_0+ t\phi_1$, $t\in [0,1]$. 

1) We have 
\begin{equation} \label{eq:Primitive}
\frac{1}{n + 1} \frac{d}{d t} E (\phi_t) =  -\int_\Omega (\phi_1-\phi_0) (dd^c \phi_t)^n, \, t \in [0,1].
\end{equation}
In particular if $u, v \in \mathcal{E}^1 (\Omega)$ and $u \leq v$, then  $0 \leq E (v) \leq E (u)$. 

Moreover we have 
\[
\frac{1}{n + 1} \frac{d^2}{d t} E (\phi_t)  = -n \int_\Omega (\phi_1-\phi_0) (dd^c ( \phi_1-\phi_0)) \wedge (dd^c \phi_t)^{n-1} \geq 0.   
\]
In particular, for any  $u, v \in \mathcal{E}^1 (\Omega)$,  the function $t \longmapsto E( (1-t) u + t v)$ is a convex function in $[0,1]$.

2) The functional $E : \mathcal{E}^1 (\Omega) \longrightarrow \R^+$ is lower semi-continuous on $\mathcal{E}^1 (\Omega)$ for the $L^1_{\rm loc}(\Omega)$-topology.
\end{lem}
The computations in 1) of the above lemma follow from integrations by part, see \cite[Corollary 3.4]{Ceg04}, and \cite{Ras17}. The second statement can be proved by following the arguments in \cite{BBGZ13}, see also \cite{ACC12}, \cite{Lu15}. 

As a consequence of the formula \eqref{eq:Primitive}, we have the following cocycle formula.
\begin{proposition} Let $u, v \in \mathcal E^1(\Omega)$. Then
$$
E(v) - E(u) = \sum_{j = 0}^n \int_\Omega (u-v) (dd^c u)^j \wedge (dd^c v)^{n-j}.
$$
\end{proposition} 

Using the above and integrations by parts we can prove the following. 
\begin{lemma} Let $u, v \in \mathcal E^1(\Omega)$. Then for any $1\leq j \leq n$
\begin{equation} \label{eq:Integration-part}
\int_{\Omega}(u - v) (dd^c u)^{n} \leq \int_{\Omega}(u - v) (dd^c u)^{j}\wedge(dd^c v)^{n - j} \leq \int_{\Omega}(u - v) (dd^c v)^{n}.
\end{equation}
As a consequence we have 
\begin{equation}\label{eq: convexity E}
 \int_{\Omega}(u - v) (dd^c u)^n \leq \frac{1}{n+1}(E(v)-E(u)) \leq \int_{\Omega} (u-v) (dd^c v)^n. 
\end{equation}
\end{lemma}

\subsection{The comparison principle}
One of the main tools in Pluripotential Theory is the maximum principle (see \cite{BT87}, \cite{BGZ09}).
\begin{proposition} \label{prop:BT1} Let $u, v \in \mathcal E^1(\Omega)$. Then
$$
{\bf 1}_{\{u> v\}} (dd^c \max \{u,v\})^n  = {\bf 1}_{\{u> v\}} (dd^c u)^n 
$$
in the sens of Borel measures on $\Omega$.
\end{proposition}
From this result we can deduce the following important consequence.

\begin{proposition} \label{prop:Dem1}  Let $u, v \in \mathcal E^1(\Omega)$. Then
$$
 (dd^c \max \{u,v\})^n  \geq {\bf 1}_{\{u>  v\}} (dd^c u)^n + {\bf 1}_{\{u\leq v\}} (dd^c v)^n
$$
in the sens of Borel measures on $\Omega$.

In particular if $v \leq u$ in $\Omega$ then
\begin{equation} \label{eq:Dem2}
{\bf 1}_{\{ v= u\}} (dd^c v)^n \leq {\bf 1}_{\{ v= u\}} (dd^c u)^n,
\end{equation}
in the sens of Borel measures on $\Omega$.
\end{proposition}

We will also need the following domination principle.  
\begin{proposition}
	Assume $u,v \in \mathcal E^1(\Omega)$ and $(dd^c u)^n \leq (dd^c v)^n$ on $\{u<v\}$. Then $u\geq v$.  
\end{proposition}
\begin{proof}
	The function $w=\max(u,v) \in \mathcal E^1(\Omega)$ satisfies $(dd^c w)^n \geq (dd^c u)^n$. Thus the comparison principle (see \cite{Ceg98}) gives $u\geq \max(u,v)\geq v$. 
\end{proof}

\subsection{Cegrell's inequalities}
The following result will be useful in the sequel.
\begin{proposition} Let $u, v, \phi \in \mathcal E^1(\Omega)$. Then

 \begin{equation}\label{eq: Cegrell inequality 2}
\int_{\Omega} (-\phi) dd^c u \wedge  (dd^c v)^{n-1} \leq \left(\int_\Omega (-\phi) (dd^c u)^n\right)^{1\slash n} \left(\int_\Omega (-\phi) (dd^c v)^n\right)^{(n - 1)\slash n}	
 \end{equation}
and 
 \begin{flalign} \label{eq: Cegrell inequality 3}
	 \int_{\Omega} (-u) (dd^c v)^n 
	\leq   E(u)^{1/(n+1)}  E(v)^{n/(n+1)}. 
\end{flalign}
\end{proposition}
The first inequality is proved in \cite{Ceg98} and the second follows from a general inequality proved in \cite[Theorem 4.6]{Per99}.
For convenience we will show how \eqref{eq: Cegrell inequality 3} can be easily deduced from an inequality of Cegrell.

\begin{proof}
		It follows from an integration by parts, see \cite[Theorem 3.2]{Ceg04}, and \cite[Theorem 5.5]{Ceg04}, that 
\begin{flalign}
	\label{eq: Cegrell inequality 1}  \nonumber  \int_{\Omega} (-u) (dd^c v)^n &= \int_{\Omega} (-v) (dd^c u)\wedge (dd^c v)^{n-1}\\
	&\leq \left ( \int_{\Omega} (-v) (dd^c u)^n \right )^{1/n}\left ( \int_{\Omega} (-v) (dd^c v)^n \right )^{(n-1)/n}. 
\end{flalign}
Changing the role of $u$ and $v$ we also get 
\begin{flalign}
	\label{eq: Cegrell inequality 2}  \nonumber  \int_{\Omega} (-v) (dd^c u)^n &= \int_{\Omega} (-u) (dd^c v)\wedge (dd^c u)^{n-1}\\
	&\leq \left ( \int_{\Omega} (-u) (dd^c v)^n \right )^{1/n}\left ( \int_{\Omega} (-u) (dd^c u)^n \right )^{(n-1)/n}. 
\end{flalign}
Plugging \eqref{eq: Cegrell inequality 2} into \eqref{eq: Cegrell inequality 1} we obtain 
\begin{flalign*}
 \int_{\Omega} (-u) (dd^c v)^n 
	\leq \left ( \int_{\Omega} (-u) (dd^c v)^n \right )^{1/n^2} E(u)^{(n-1)/n^2} E(v)^{(n-1)/n},
\end{flalign*}
which yields the result for test functions. The general case follows by approximation.
\end{proof}

\subsection{Plurisubharmonic envelopes}
Plurisubharmonic envelopes are classical objects and play an important role in (Pluri)-potential Theory. They were used in  \cite{BT76} to solve complex Monge-Amp\`ere equations via the Perron method.

More recently they were studied  in \cite{GLZ19} and applied to solve some highly  degenerate complex Monge-Amp\`ere  equations. Let us recall the definition and the main properties following the presentation in \cite{GLZ19}.

Let $h : \Omega \longrightarrow \R$ be a   Borel function bounded from above in $\Omega$. We define its plurisubharmonic envelope in $\Omega$ by the formula
$$
P(h) := \left(\sup \{v \in \PSH(\Omega) \, ; \, v \leq h \, \, \text{in} \, \, \Omega\}\right)^*,
$$
provided there exist $v_0 \in \PSH(\Omega)$ such that $v_0 \leq h$ in $\Omega$.
We have the following important result (\cite{GLZ19}).
\begin{theorem} \label{thm:projection} Let $h : \Omega \longrightarrow \R^-$ be a negative quasi-continuous function.
Assume there exists $v_0 \in \mathcal E^1(\Omega)$ such that $v_0 \leq h$ in $\Omega$.
Then $u := P(h) \in \mathcal E^1(\Omega)$,  $u \leq h$ quasi everywhere in $\Omega$ and $(dd^c u)^n$ is concentrated on the contact set
$\{ z \in \Omega ; u(z) = h(z)\}$.
\end{theorem}
This theorem was proved in the compact K\"ahler setting in \cite{GLZ19}. The proof of the above theorem is done in the same way, see \cite{ACLR25}. 

\section{The eigenvalue problem} \label{sect: eigenvalue problem}
Let $\Omega$  be a bounded hyperconvex domain in $\mathbb C^n$, $\mu$ be a positive (locally finite) Borel measure on $\Omega$ with $\mu(\Omega) > 0$, and let  $I_{\mu}$ denote the functional  $u\mapsto  \int_{\Omega}|u|^{n+1} d\mu$ with values in $[0,+\infty]$. If $u\in \mathcal E^1(\Omega) \setminus \{0\}$, then $I_{\mu}(u)>0$, but in general it may take the value $+\infty$. 

The complex Monge-Amp\`ere eigenvalue problem for the datum $(\Omega,\mu)$ consists of finding all couples $(\lambda,u)\in (0,+\infty)\times \mathcal E^1(\Omega)\setminus \{0\}$ such that
\begin{equation*} \label{eq: MA eigenvalue}\tag{${\rm MA}_{\mu,\lambda}$}
 (dd^c u)^n = (-\lambda u)^n \mu.
\end{equation*}
If $\lambda>0$ and \eqref{eq: MA eigenvalue} has a non-trivial solution $u\in \mathcal{E}^1(\Omega)$, then $\lambda$ is called an eigenvalue, and $u$ is called an eigenfunction. 

A couple $(\lambda,u)\in (0,+\infty) \times \mathcal E^1(\Omega) \setminus \{0\}$ is a subsolution of \eqref{eq: MA eigenvalue} if 
\[
(dd^c u)^n \geq (-\lambda u)^n \mu,
\] 
in the sense of Radon measures. In this case we also say that $u$ is a subsolution to \eqref{eq: MA eigenvalue}.  It is called a supersolution of \eqref{eq: MA eigenvalue} if the reverse inequality holds.

If $(\lambda,u)$ is a subsolution of \eqref{eq: MA eigenvalue} then $I_{\mu}(u)$ must be finite since 
\[
\lambda^n I_{\mu}(u) = \int_{\Omega} (-u) (-\lambda u)^n d\mu \leq \int_{\Omega}  (-u) (dd^c u)^n =E(u)<+\infty. 
\]

\subsection{Rayleigh quotient type formula}

\begin{definition}
If there exists at least one $v\in \mathcal{E}^1(\Omega)$ such that $I_{\mu}(v)<+\infty$, then we define
\[
 	\lambda_1(\mu) := \inf \left \{ \frac{E(u)}{I_{\mu}(u)} \; , \; u \in \mathcal E^1(\Omega) \setminus \{0\}, \; I_{\mu}(u)<+\infty \right \}. 
 \]
\end{definition}

Observe that $I_\mu (u) > 0$ if $ u \in \mathcal E^1(\Omega)\setminus \{0\}$. Indeed take a compact set $K \subset \Omega$ such that $\mu(K) >0$ and set $m := \max_K u < 0$.
Then 
$$
I_\mu (u) \geq (-m)^{n+1} \mu(K) >0.
$$
 
In the sequel we will simply write $\lambda_1$ if the measure $\mu$ is fixed. 

\begin{lemma}\label{lem: def of lambda 1}
 Assume $\varphi\in \mathcal E^1(\Omega)\setminus\{0\}$ and $I_{\mu}(\varphi)<+\infty$. Then 
 \[
 \lambda_1(\mu)=  \lambda_1(\mu,\varphi) := \inf \left \{ \frac{E(u)}{I_{\mu}(u)} \; , \; u \in \mathcal E^1(\Omega), \; \varphi \leq u <0 \right \}. 
 \] 
 In particular, if $I_{\mu}$ is finite on $\mathcal E^1(\Omega)$, then 
  \[
 \lambda_1(\mu)= \inf \left \{ \frac{E(u)}{I_{\mu}(u)} \; , \; u \in \mathcal E^1(\Omega) \setminus \{0\} \right \}. 
 \] 
\end{lemma}
\begin{proof}
Fix $u\in \mathcal E^1(\Omega)$ with $\varphi\leq u <0$. Then  $I_{\mu}(u)<+\infty$, thus $\lambda_1(\mu) \leq \lambda_1(\mu,\varphi)$. 

	For $t>0$, consider $u_t:= \max(u,t\varphi)$. Then $u_t\in \mathcal E^1(\Omega)\setminus \{0\}$ and $u_t\geq t\varphi$. Since $\lambda_1(\mu,\varphi)= \lambda_1(\mu,t\varphi)$, we also have 
	\[
	\frac{E(u_t)}{I_{\mu}(u_t)} \geq \lambda_1(\mu,\varphi)^n.
	\]
	Letting $t\to +\infty$, and using the continuity of $E$ as well as the monotone convergence theorem we arrive at 
	\[
	\frac{E(u)}{I_{\mu}(u)} \geq \lambda_1(\mu,\varphi)^n.
	\]
	Therefore $\lambda_1(\mu) \geq \lambda_1(\mu,\varphi)$, and the equality holds.  
\end{proof}

\begin{lemma}\label{lem: envelope is supersolution}
	Assume $u$ is a supersolution, while $v$ is a subsolution of \eqref{eq: MA eigenvalue}.
Then $\varphi := P(u-v,0)$ is a supersolution of \eqref{eq: MA eigenvalue}.
\end{lemma}
Here, and in the sequel, we use the notation $P(f,g):=P(\min(f,g))$.
\begin{proof}
	By \cite{GLZ19} since $u \leq \varphi \leq 0$ in $\Omega$ we know that $\varphi \in \mathcal E^1(\Omega)$ and $(dd^c \varphi)^n$ is concentrated on the contact set 
$\mathcal C :=\{\varphi=u-v\}$.  Since $\varphi + v \leq u$ in $\Omega$, it follows from the comparison principle that 
$$ 
(dd^c \varphi)^n = {\bf 1}_{ \mathcal C } (dd^c \varphi)^n  \leq {\bf 1}_{ \mathcal C }(dd^c (\varphi + v))^n \leq {\bf 1}_{ \mathcal C }(dd^c u)^n.
$$ Hence $(dd^c \varphi)^n$ is absolutely continuous with respect to $(dd^c u)^n\leq (-\lambda u)^n \mu$. We can thus write $(dd^c \varphi)^n = h^n \mu$ where $0 \leq h$ is such that $h^n \in L^1_{\rm loc}(\Omega,\mu)$. Then the mixed Monge-Amp\`ere inequalities, see \cite{Dinew09Z}, give, since $\varphi+v \leq u$,
\[
{\bf 1}_{\{\varphi=u-v\}}(dd^c u)^n \geq {\bf 1}_{\{\varphi=u-v\}}(dd^c v+dd^c \varphi)^n\geq {\bf 1}_{\{\varphi=u-v\}}(-\lambda v + h)^n \mu.  
\]
It follows that 
\[
{\bf 1}_{\{\varphi=u-v\}}(-\lambda u)^n d\mu \geq {\bf 1}_{\{\varphi=u-v\}}(-\lambda v + h)^n d\mu.
\]
We thus have
\[
{\bf 1}_{\{\varphi=u-v\}}  h \leq  {\bf 1}_{\{\varphi=u-v\}} \lambda (v-u),
\]
almost everywhere with respect to $\mu$. 
Therefore, using that $(dd^c \varphi)^n$ vanishes outside $\{\varphi=u-v\}$, we have
\[
(dd^c \varphi)^n \leq {\bf 1}_{\{\varphi=u-v\}} (\lambda (v-u))^n \mu \leq (-\lambda  \varphi)^n \mu. 
\]
Therefore, $\varphi\in \mathcal E^1(\Omega)$ is a supersolution to \eqref{eq: MA eigenvalue}.  
\end{proof}

From the definition of $\lambda_1$ we deduce the following simple but very useful observation.

\begin{lemma}\label{lem: supersol is sol}
	Assume $u \in \mathcal E^1(\Omega)\setminus \{0\}$ and $I_{\mu}(u)<+\infty$.  If $(\lambda_1,u)$ is a supersolution of \eqref{eq: MA eigenvalue}, then $(\lambda_1,u)$ is a solution to \eqref{eq: MA eigenvalue}. 
\end{lemma}
\begin{proof} 
Since $(\lambda_1,u)$ is a supersolution, we  have 
\[
(-u)(dd^c u)^n \leq \lambda_1^n (-u)^{n+1} \mu.
\] 
Integrating over $\Omega$, and using the definition of $\lambda_1$, we thus get 
\begin{flalign*}
	E(u) = \int_{\Omega} (-u) (dd^c u)^n \leq \lambda_1^n\int_{\Omega} (-u)^{n+1} d\mu \leq E(u).
\end{flalign*} 
	All inequalities above are thus equalities, in particular we have 
	\[
	(-u)(dd^c u)^n = \lambda_1^n(- u)^{n+1} \mu.
	\]
	Since $u<0$, we can divide by $(-u)$ to get $(dd^c u)^n = \lambda_1^n(- u)^n \mu$, finishing the proof.  	 
\end{proof}

\begin{theorem}\label{thm: uniqueness of the eigenvalue}
	Assume $\varphi \in \mathcal E^1(\Omega) \setminus \{0\}$ and $I_{\mu}(\varphi)<+\infty$. If $\varphi$ is a subsolution to \eqref{eq: MA eigenvalue}, then $\lambda \leq \lambda_1$. If $\varphi$ solves \eqref{eq: MA eigenvalue} then $\lambda=\lambda_1$.  
\end{theorem}
\begin{proof}
Fix a closed ball $K\subset \Omega$, a constant $\varepsilon>0$ such that $\varphi < -2\varepsilon$ on $K$. Let $\mathcal S(\varphi,K,\varepsilon)$ denote the set of all $u\in \PSH(\Omega)$ such that $\varphi \leq u < 0$ on $\Omega$ and  $u \leq -\varepsilon$ on $K$.  Define 
\[
\gamma_1^n:=  \inf \left\{ \frac{E(u)}{I_{\mu}(u)} \; : \; u\in \mathcal S(\varphi,K,\varepsilon) \right \}.
\]
Let $(u_j)$ be a sequence of psh functions in $\Omega$, such that $\varphi \leq u_j<0$ on $\Omega$, $u_j\leq -\varepsilon$ on $K$, and 
\[
\lim_{j\to +\infty} \frac{E(u_j)}{I_{\mu}(u_j)} = \gamma_1^n.
\]
 Extracting a subsequence we can assume that $u_j \to u$ in $\Omega$ almost everywhere with respect to $dV$ and in $L^1(dV)$. The condition $u_j\leq -\varepsilon$ on $K$ yields $u<0$. By semicontinuity of $E$ we have 
 \[
 \liminf_{j\to +\infty} E(u_j) \geq E(u).
 \]
  Since $u_j\geq \phi$, by the dominated convergence theorem we also have $I_{\mu}(u_j) \to I_{\mu}(u)$, hence 
  \[
  E(u)- \gamma_1^n I_{\mu}(u)=0 \leq E(v) -\mu_1 I_{\mu}(v),
  \]
   for all $v \in \mathcal S(\varphi,K,\varepsilon)$.  We now claim that 
\[
 (dd^c u)^n \leq \gamma_1^n (-u)^n \mu \; \text{in} \; U:= \Omega \setminus K \cup \{u<-\varepsilon\}.
 \]    
 Fix a positive test function $\chi$ which is compactly supported in $U$. Then,  for sufficiently small $t>0$, we have $u_t:= P(u+t\chi) \in \mathcal S(\varphi,K,\varepsilon)$. The function 
   \[
	g(t):= E(u_t) - \gamma_1^n I_{\mu}(u+t\chi) 
   \]
   therefore satisfies $g(t) \geq g(0)$ for all $t\in [0,\delta)$, where $\delta$ is a small constant. It thus follows that 
   $g'(0^+)\geq 0$, and this can be rewritten as 
   \[
   \lim_{t\to 0^+} \frac{E(u_t)-E(u)}{t} \geq -(n+1) \gamma_1^n \int_{\Omega} \chi (-u)^n d\mu. 
   \]
   Using 
   \[
   \frac{E(u_t)-E(u)}{t} \leq (n+1)\int_{\Omega} \frac{u-u_t}{t} (dd^c  u_t)^n = (n+1)\int_{\Omega} (-\chi) (dd^c  u_t)^n,
   \]
   we then arrive at 
   \[
   \int_{\Omega} \chi (dd^c u)^n \leq \int_{\Omega} \chi (-\gamma_1 u)^n d\mu.
   \]
   Since this is true for all such $\chi$, the claim follows. 
   
We next prove that $\gamma_1 \geq \lambda$. Assume by contradiction that it is not the case. 
 Then $a= \min(\lambda\gamma_1^{-1},2)>1$, and we have the inclusion $\{au<\varphi\} \subset U$. Therefore, in $\{au<\varphi\}$, we have
 \[
 (dd^c a u)^n  \leq (-a\gamma_1 u)^n \mu \leq  (-\lambda u)^n \mu  \leq (-\lambda \varphi)^n \mu \leq  (dd^c \varphi)^n.
 \]
 By the domination principle we thus have $au \geq \varphi$ in $\Omega$. It follows that $au$ minimizes $E - \gamma_1^n I_{\mu}$ over all candidates in $\mathcal S(\varphi,K,\varepsilon)$. Since $au<- \varepsilon$ on $K$, we have $\Omega\setminus K \cup \{au<-\varepsilon\}=\Omega$, and the previous step yields 
 \[
 (dd^c au )^n \leq (-\gamma_1 au )^n \mu\; \text{in} \; \Omega. 
 \]
 By the definition of $\gamma_1$, we must have $(dd^c u )^n = (-\gamma_1 u )^n \mu$. 
 
 We repeat the same argument with $a^2 u$ and $\varphi$ to see that on $\{a^2u <\varphi\}$, we have 
 \[
 (dd^c a^2 u)^n = (-a^2 \gamma_1 u)^n \mu \leq  (-\lambda a u)^n \mu \leq (-\lambda \varphi)^n \mu =(dd^c \varphi)^n.
 \]
 The domination principle then gives $a^2u \geq \varphi$. Repeating the argument we see that $a^ju\geq \varphi$ for all $j\in \mathbb N$, which is impossible. We therefore have $\gamma_1\geq \lambda$. Letting $\varepsilon\to 0$ we obtain $\lambda_1\geq \lambda$.

Finally, we assume that $(\lambda,\varphi)$ is a solution to \eqref{eq: MA eigenvalue}. Since $\varphi$ is a candidate defining $\lambda_1(\mu,\varphi)$, we also have $\lambda \geq \lambda_1(\mu,\varphi)$, therefore $\lambda=\lambda_1(\mu,\varphi)$. 
\end{proof}

As a consequence of the uniqueness of the eigenvalue, we obtain a necessary condition for the solvability of \eqref{eq: MA eigenvalue}. 
\begin{proposition}
	\label{prop: finiteness}
	If the eigenvalue problem \eqref{eq: MA eigenvalue} has a solution then we have $\mathcal E^1(\Omega) \subset L^{n+1}(\Omega,d\mu)$. In particular, 
  \[
 \lambda_1(\mu)= \inf \left \{ \frac{E(u)}{I_{\mu}(u)} \; , \; u \in \mathcal E^1(\Omega) \setminus \{0\} \right \}. 
 \] 
\end{proposition}
\begin{proof}
	Let $(\lambda_1,u)$ be a solution to \eqref{eq: MA eigenvalue}, with $\lambda_1>0$, $u\in \mathcal E^1(\Omega)$, $u<0$. If $v\in \mathcal E^1(\Omega)$ and $t>0$, then $\max(v,tu)\in \mathcal E^1$ and $I_{\mu}(\max(v,tu)<+\infty$. Theorem \ref{thm: uniqueness of the eigenvalue} thus yields
	\[
	\lambda_1^n \int_{\Omega} |\max(v,tu)|^{n+1} d\mu \leq E(\max(v,tu)) \leq E(v). 
	\]
	Letting $t\to +\infty$ and using the monotone convergence theorem, we arrive at the result. 
\end{proof}

We end this subsection with the following inequality which may be of independent interest. 
\begin{corollary}\label{cor: an inequality in E1}
	If $u,v\in \mathcal E^1(\Omega)\setminus \{0\}$, then 
	\[
	\int_{\Omega} (-u)^{-n}(-v)^{n+1} (dd^c u)^n \leq E(v). 
	\]	
\end{corollary}
\begin{proof}
	We consider $\mu= (-u)^{-n}(dd^c u)^n$. Then $1$ is an eigenvalue and $u$ is an eigenfunction of the eigenvalue problem \eqref{eq: MA eigenvalue}. By Theorem \ref{thm: uniqueness of the eigenvalue} and Proposition \ref{prop: finiteness}, we have $I_{\mu}(v) \leq E(v)$, finishing the proof.
\end{proof}

\subsection{Uniqueness of eigenfunctions}

\begin{theorem}\label{thm: subsol is sol}
	Assume $\lambda>0$, $\varphi,\psi\in \mathcal E^1\setminus \{0\}$. If $\varphi$ is a solution and $\psi$ is a subsolution to \eqref{eq: MA eigenvalue}, then $\psi$ is also a solution and $\varphi=c\psi$ for some constant $c>0$.
\end{theorem}
\begin{proof}
From Theorem \ref{thm: uniqueness of the eigenvalue}  we have $\lambda=\lambda_1(\mu,\varphi)$. 

We first assume that $(dd^c \psi)^n$ is absolutely continuous with respect to $\mu$. Let $a>0$ be a constant such that $\{\varphi-a\psi<0\}$ is not empty, and set $v:=P(\varphi-a\psi,0)$. Then $v<0$ in $\Omega$ since $v$ is plurisubharmonic in $\Omega$, $v\leq 0$, and it is not identically $0$. It follows from Lemma \ref{lem: envelope is supersolution} that $(\lambda, v)$ is a supersolution to \eqref{eq: MA eigenvalue}. Since $v\geq \varphi$, it follows from Lemma \ref{lem: supersol is sol} that $(\lambda,v)$ solves \eqref{eq: MA eigenvalue}. Therefore, using the fact that $(dd^c v)^n$ is supported on the set $\{v=\varphi-a\psi\}$, which is contained in the set $\{\varphi<a\psi\}$, we have 
	\[
	 \int_{\{a\psi<\varphi\}} (dd^c v)^n =0.
	\]
	This implies $\mu(a\psi<\varphi)=0$ since $v<0$ and $(dd^c v)^n =(-\lambda v)^n \mu$. On the other hand, since $(dd^c \psi)^n$ is absolutely continuous with respect to $\mu$, we also have $(dd^c a\psi)^n (a\psi<\varphi) =0$, and the domination principle yields $a\psi\geq \varphi$.  
	
	Now let $c>0$ be the supremum of all $a>0$ such that $\{\varphi-a\psi<0\}$ is not empty. By the previous argument we have $\varphi \leq a\psi$ for all such $a$, hence $\varphi \leq c \psi$. From this we infer that $c$ is finite. From the definition of $c$ we also have $\varphi \geq  c \psi$, therefore  $\varphi =c \psi$. The latter implies that $(\lambda,\psi)$ solves \eqref{eq: MA eigenvalue}, finishing the first step.

	To treat the general case, observe that by H\"older inequality, the positive Borel measure $\nu :=  (-\lambda \psi)^n \mu$ satisfies the following a priori estimate : for any $\phi \in \mathcal E^1(\Omega)$ we have
$$
\int_\Omega (-\phi) d \nu \leq A \left(\int_\Omega (-\phi)^{n+1} d \mu\right)^{1\slash n+1}, \, \, \text{with} \, \, A := \lambda^n \left(\int_\Omega (-\psi)^{n+1} d \mu\right)^{n\slash n+1}\cdot
$$
 By \cite[Theorem A]{ACC12}, there exists $u\in \mathcal E^1(\Omega)$ such that $(dd^c u)^n = (-\lambda \psi)^n \mu$. Since $(dd^c \psi)^n \geq (-\lambda \psi)^n \mu \geq (dd^c u)^n$, the comparison principle ensures that $u\geq \psi$, hence $(\lambda,u)$ is also a subsolution to \eqref{eq: MA eigenvalue}.  It thus follows from Lemma \ref{lem: supersol is sol} that $(\lambda,u)$ solves \eqref{eq: MA eigenvalue}. Now, assume by contradiction that $(\lambda,\psi)$ is not a solution to \eqref{eq: MA eigenvalue}.  Then we can write $(dd^c \psi)^n = (-\lambda \psi)^n (\mu + \sigma)$, where the measure $\sigma$ is positive and non-pluripolar.  From Theorem \ref{thm: uniqueness of the eigenvalue}  we have 
	\[
	\lambda = \lambda_1(\mu,\varphi) = \lambda_1(\mu+\sigma,\psi). 
	\]
Using that $0>u\geq \psi$, we now arrive at 
\[
\lambda^n = \frac{E(u)}{I_{\mu}(u)} > \frac{E(u)}{I_{\mu+\sigma}(u)} \geq (\lambda_1(\mu+\sigma))^n= \lambda^n,
\]
which is a contradiction. 

Therefore, $(\lambda,\psi)$ also solves \eqref{eq: MA eigenvalue}, and by the first step, $\varphi=c\psi$ for some positive constant $c$. 
\end{proof}

We now prove the uniqueness.

\begin{corollary}\label{cor: unique}
	Assume $\lambda>0$, $\gamma>0$, and $\varphi,\psi\in \mathcal E^1\setminus \{0\}$. If $(\lambda,\varphi)$  and $(\gamma,\psi)$ satisfy the following equation
$$
 (dd^c \varphi)^n = (-\lambda \varphi)^n \mu, \, \, \, \, (dd^c \psi)^n = (-\gamma \psi)^n \mu
$$
then  $\lambda=\gamma$ and $\varphi=c\psi$, for some positive constant $c$.  
\end{corollary}

\begin{proof}
	Since $\lambda$ and $\gamma$ play the same role, we can assume that $\lambda\leq \gamma$. By Theorem \ref{thm: uniqueness of the eigenvalue} we have $\lambda=\lambda_1(\mu,\varphi)$, and by Theorem \ref{thm: subsol is sol}, we know that $(\lambda,\psi)$ also solves \eqref{eq: MA eigenvalue}. It thus follows that $\lambda = \gamma$ and $\varphi=c\psi$, for some positive constant $c$.
\end{proof}

\subsection{Existence of eigenfunctions}
The goal of this section is to give a sufficient condition for $\mu$ which ensures that \eqref{eq: MA eigenvalue} has a solution. As we have seen in the previous section, there is at most one eigenvalue and all the eigenfunctions are proportional.

\begin{lemma}\label{lem: quantitative bound}
	Assume that $I_{\mu}$ is finite on $\mathcal E^1$. Then there exists a constant $A>0$ such that, for all $u\in \mathcal E^1$, $I_{\mu}(u)\leq A\; E(u)$.  
\end{lemma}
We recall that $E(u):= \int_{\Omega} (-u) (dd^c u)^n$, for $u\in \mathcal E^1$. 
\begin{proof}
Assume by contradiction that there exists a sequence $(u_j)$ in $\mathcal E^1$ such that $I_{\mu}(u_j)\geq 4^j E(u_j)$. Define 
\[
v_{k} = \sum_{j =1}^k \frac{u_j}{2^j E(u_j)^{1/(n+1)}}. 
\]
By the convexity of $E$, we have that $E(v_k) \leq 1$. Thus $v_k$ decreases to $v\in \mathcal E^1$. By the monotone convergence theorem, we also have $I_{\mu}(v) =+\infty$, contradicting the finiteness assumption on $I_{\mu}$. 
\end{proof}

\begin{example}
	If $\mu = (dd^c \phi)^n$ for some bounded plurisubharmonic function $\phi$, then by B{\l}ocki's inequality \cite{Bl93}, $I_{\mu}$ is finite on $\mathcal E^1$. 
\end{example}

We let $\mathcal E^1_C$ denote the set of functions $u\in \mathcal E^1$ with $E(u)\leq C$. We recall that $\mathcal E^1_C$ equipped with the $L^1$-topology is compact. 

 The variational method of \cite{BBGZ13} can be applied exactly as in \cite{BZ23} to solve the equation \eqref{eq: MA eigenvalue}. 
\begin{theorem}\label{thm: existence} \cite{BZ23}
	Assume that $I_{\mu}$ is continuous on $\mathcal E^1_C$, for some $C>0$. Then the equation \eqref{eq: MA eigenvalue} has a solution $(\lambda_1,u)$. 
\end{theorem} 
In Theorem \ref{thm: existence iterative} we provide a different proof inspired by \cite{Zer25}.

\begin{example}
	Assume $\mu = f dV$, with $\int_{\Omega} f (\log f)^n (\log \log f)^p dV<+\infty$, for some $p>0$. Then $\mu$ satisfies the assumption of Theorem \ref{thm: existence}. Indeed, let $(u_j)$ be a sequence in $\mathcal E^1_C$ converging in $L^1(dV)$ to $u\in \mathcal E^1_C$. By \cite{DGL21}, there exists $\gamma>0$ such that 
	\[
	\int_{\Omega} e^{-\gamma |u_j|^{1+1/n}} dV \leq A. 
	\] 
	Given $B>1$, it follows that 
	\begin{flalign*}
		\int_{\{f>B\}} |u_j|^{n+1} fdV  &= \int_{\{|u_j|^{1+1/n}< \log f\} \cap \{ f>B\}} f (\log f)^n dV\\
		&+  A_1 \int_{\{|u_j|^{1+1/n} \geq \log f\} \cap \{ f>B\}}  e^{2^{-1}\gamma |u_j|^{1+1/n}} |u_j|^{n+1} dV \\
		& = o(B). 
	\end{flalign*}
	From this we deduce that $\int_{\Omega} |u_j|^{n+1} fdV\to \int_{\Omega} |u|^{n+1} fdV$.

	 When $p\leq n-1$ the solution to $(dd^c \phi)^n= fdV$ may fail to be bounded, see \cite{GL25}. Therefore,  in this case $\mu=fdV$ does not satisfies the assumptions in \cite{BZ23}, \cite{BZ24}. 
\end{example}

The continuity of $I_{\mu}$ implies in fact a stronger convergence.

\begin{lemma}\label{lem: convergence}
	Assume that $I_{\mu}$ is continuous on $\mathcal E^1_C$, for some $C>0$, and let $(u_j)$, $(v_j)$ be sequences in $\mathcal E^1_C$ converging in $L^1(dV)$ to $u, v\in \mathcal E^1_C$. Then 
	\[
	\int_{\Omega} ||u_j|^{n+1}-|u|^{n+1}|d\mu \to 0. 
	\] 
	and, for all $0\leq p\leq n+1$,
	 \begin{equation}
	  	\label{eq: convergence I mu}
	  	  \lim_{j\to +\infty} \int_{\Omega} (-u_j)^p (-v_j)^{n+1-p} d\mu = \int_{\Omega} (-u)^p (-v)^{n+1-p}d\mu.  
	  \end{equation}
\end{lemma}

\begin{proof}
	 For each $j$, we define $\tilde{u}_j=(\sup_{k\geq j} u_k)^*$. Then $\tilde{u}_j\searrow u$ pointwise in $\Omega$. By classical properties of plurisubharmonic functions we have $\limsup_j u_j \leq u$ in $\Omega$. For all $p,q\geq 0$, with $p+q\leq n+1$, Fatou's lemma thus yields
	 \[
	 \liminf_{j\to +\infty} \int_{\Omega} (-u_j)^p (-\tilde{u}_j)^q (-u)^{n+1-p-q} d\mu \geq \int_{\Omega} (-u)^{n+1}d\mu. 
	 \]
	 Summing up the above inequalities (after multiplying with appropriate binomial coefficients) we obtain 
	 \[
	 \liminf_{j\to +\infty} \int_{\Omega}  (-u_j -\tilde{u}_j- u)^{n+1}d\mu \geq \int_{\Omega} (-3u)^{n+1}d\mu. 
	 \]
	 Since $I_\mu$ is continuous in $\mathcal E^1_C$ (hence also in $\mathcal E^1_{3C}$ by scaling invariant), we deduce that the above inequality is in fact an equality. We thus infer that 
	  \begin{equation}
	  	\label{eq: convergence proof 1}
	  	  \lim_{j\to +\infty} \int_{\Omega} (-u_j)^p (-\tilde{u}_j)^q (-u)^{n+1-p-q} d\mu = \int_{\Omega} (-u)^{n+1}d\mu.  
	  \end{equation}
	 Since $\tilde{u}_j\geq \max(u_j,u)$, we can write 
	 \[
	 |u_j-u| = 2\max(u_j,u)- u_j-u \leq 2(\tilde{u}_j -u) + u-u_j,
	 \]
	 and hence
	 \begin{flalign*}
	 	 \int_{\Omega} ||u_j|^{n+1}-|u|^{n+1}|d\mu &= \sum_{p=0}^n \int_{\Omega} |u_j-u| (-u_j)^p (-u)^{n-p} d\mu \\
	 	 &\leq  \sum_{p=0}^n \int_{\Omega} (2(\tilde{u}_j-u)+u-u_j) (-u_j)^p (-u)^{n-p} d\mu.
	 \end{flalign*}
	 Using \eqref{eq: convergence proof 1}, we see that the right-hand side tends to $0$ as $j\to +\infty$, finishing the proof of the first statement. To prove \eqref{eq: convergence I mu} we note that, by Fatou's lemma
	 \[
	 \liminf_j \int_{\Omega} (-u_j)^p (-v_j)^{n+1-p} d\mu \geq \int_{\Omega} (-u)^p (-v)^{n+1-p} d\mu, \forall p\leq n+1.
	 \]
	 On the other hand, since $I_{\mu}$ is continuous on $\mathcal E^1_C$, for all $C>0$ fixed, we have 
	 \[
	 \lim_j \int_{\Omega} (-u_j+v_j)^{n+1} d\mu = \int_{\Omega} (-u-v)^{n+1}d\mu. 
	 \]
	 It thus follows that all the above inequalities are equalities, finishing the proof. 
\end{proof}

As an immediate consequence we have the following
\begin{corollary}\label{cor: I cont}
		Fix $C>0$. If $I_\sigma$ is continuous in $\mathcal E^1_C$ and $\mu \leq \sigma$, then $I_{\mu}$ is also continuous in $\mathcal E^1_C$. 
\end{corollary}

We now give a sufficient condition for the continuity of $I_{\mu}$ on $\mathcal E^1_C$. 

\begin{prop} \label{prop: I cont}
	If $\mu\leq  (dd^c \phi)^n$, for some $\phi\in \PSH(\Omega)\cap C^0(\overline\Omega)$, then $I_{\mu}$ is continuous on $\mathcal E^1_C$ for any $C>0$. 
\end{prop}
\begin{proof}
	Fix $C>0$. By Corollary \ref{cor: I cont}, we can assume that $\mu= (dd^c \phi)^n$, for some $\phi\in \PSH(\Omega)\cap C^0(\overline \Omega)$. Let $(f_k)$ be a sequence of smooth functions decreasing to $\phi$, and set $\phi_k = P(f_k,0)$. Then  $\phi_k\in \PSH(\Omega)\cap C^0(\overline \Omega)$,  $\phi_k|_{\partial \Omega}=0$,  $\phi_k \searrow \phi$ and $\mu_k:=(dd^c \phi_k)^n\leq A_k dV$. The latter implies that, for each $k$, $I_{\mu_k}$ is continuous in $\mathcal E^1_C$.  It is thus sufficient to prove that $I_{\mu_k}$ uniformly converges to $I_{\mu}$. To do this, we first prove the following   
	\begin{equation}
		\label{eq: I cont}
		I_{\mu}(u-v) - I_{\nu}(u-v) \leq B \sup_{\Omega}|\varphi-\psi| (E(u)+ E(v)), \; \forall u, v\in \mathcal E^1_C, \; u\leq v, 
	\end{equation}
	where $\mu=(dd^c \varphi)$, $\nu=(dd^c \psi)$, $\varphi,\psi \in \PSH(\Omega)\cap L^{\infty}$, and  $B$ is a uniform constant. 
	
	By approximation, we can assume that $u,v\in \mathcal E^0$. 
	By replacing $v$ with $\max(u,v-\varepsilon)$ we can further assume $w:=u-v=0$ near $\partial \Omega$.  Then 
	\begin{flalign*}
		I_{\mu}(w) &- I_{\nu}(w) = \int_{\Omega} (-w)^{n+1} dd^c (\varphi-\psi) \wedge T,\\
		& = -(n+1)\int_{\Omega} (\varphi-\psi) (-w)^n dd^c w \wedge T \\
		&+ n(n+1) \int_{\Omega} (\varphi-\psi) (-w)^{n-1} dw \wedge d^c w \wedge T,\\
		&\leq B_1\sup_{\Omega}|\varphi-\psi| \int_{\Omega} \left ((-w)^n dd^c (u+v) 
		+ (-w)^{n-1} dw \wedge d^c w \right) \wedge T\\
		& \leq B_2 \sup_{\Omega}|\varphi-\psi| \int_{\Omega} (-w)^n dd^c (u+v)  \wedge T.
	\end{flalign*}
	where $T= \sum_{p=0}^{n-1} (dd^c \varphi)^p \wedge (dd^c \psi)^{n-1-p}$, and in the last line we integrate by parts 
	\[
	\int_{\Omega} d(-w)^n \wedge d^c w \wedge T = -\int_{\Omega} (-w)^n dd^c w  \wedge T \leq \int_{\Omega} (-w)^n dd^c (u+v)  \wedge T.
	\]
	Here $B_1,B_2$ are uniform constants. 
	 Integrating by parts several times as in \cite{Bl93}, we obtain 
	 \[
	 \int_{\Omega} (-w)^n dd^c (u+v) \wedge T \leq B_3 (\sup_{\Omega} |\varphi+\psi|)^{n-1}E(u+v). 
	 \]
	 Finally we observe that, by Dini's theorem, $\phi_k$ uniformly converges to $\phi$ in $\overline \Omega$. The estimate \eqref{eq: I cont} with $v=0$ thus completes the proof. 
\end{proof}

\section{A general Dirichlet Problem}\label{sect: Dirichlet problem}

\subsection{Monge-Amp\`ere equation with zero boundray values}
In this section we assume $\mu(\Omega)<+\infty$,  $\mathcal E^1(\Omega) \subset L^{n+1}(\mu)$, and we consider the following complex Monge-Amp\`ere  equation:
\begin{equation} \label{eq:DP}
	(dd^c u)^n  -  F(\cdot,u)^n \mu = 0, \; u \in \mathcal E^1(\Omega), 
\end{equation}
where $F :  \Omega \times (-\infty,0] \longrightarrow \R^+$ is a Borel function such that 
\begin{itemize}
\item $ f := F(\cdot,0) \in L^{n+1}(\Omega,\mu)$, and 
\item for $\mu$-a.e. $z\in \Omega$, $t \longmapsto F(z,t)$ is continuous in $]0,\infty[$. 
\end{itemize}

\begin{theorem}  \label{thm:Existence-Uniqueness} Assume that there exists $ \lambda_0 > 0 $ such that
$$
\frac{\partial F}{\partial t} \geq - \lambda_0 , \, \, \mathrm{in} \, \,  ]-\infty, 0] \times \Omega.
$$ 
Then  

(1) if $\lambda_0 < \lambda_1$,  the maximum principle holds for the operator 
$(dd^c u)^n - F(\cdot,u) \mu$ i.e. if $u \in \mathcal E^1(\Omega)$ is a supersolution and $v \in \mathcal E^1(\Omega)$ is a subsolution then $u \geq v$ in $\Omega$. 

Moreover the Dirichlet problem \eqref{eq:DP} admits a unique solution $\varphi \in \mathcal E^1(\Omega)$. In particular, if $F(z,0)\equiv 0$, the unique solution is identically $0$.

(2) If $\lambda_0 = \lambda_1$ and $u, v \in \mathcal E^1(\Omega)$ are two solutions to \eqref{eq:DP}  such that $u(z_0) < v(z_0)$ for some point $z_0 \in \Omega$, then $\varphi := P(u-v,0)<0$ is an eigenfunction associated  to the eigenvalue $\lambda_1$.

\end{theorem}
\begin{proof} 

1. We first prove the maximum principle.  Let $u, v \in \mathcal E^1(\Omega)$ such that
$(dd^c u)^n \leq F(\cdot,u) \mu$ and  $(dd^c v)^n \geq F(\cdot,v) \mu$ on $\Omega$.

To compare $u$ and $v$ we will consider $u-v$ and set $\varphi := P(u-v,0)$.
Since $u \leq u - v$ it follows that $u \leq \varphi \leq 0$ hence $\varphi \in \mathcal E^1(\Omega)$.

We claim that $\varphi$ is a supersolution to the equation \eqref{eq:DP}.

Indeed $(dd^c \varphi)^n$ is concentrated in the contact set $D := \{\varphi =u - v\} = \{\varphi + v =u\}$.
Since $\varphi \leq u - v$ we have $\varphi + v \leq u$ in $\Omega$.
It follows from the comparison principle that
\begin{eqnarray*}
(dd^c \varphi)^n = {\bf 1}_D (dd^c \varphi)^n &\leq & {\bf 1}_D (dd^c (\varphi + v))^n  \\
& \leq & {\bf 1}_D (dd^c u)^n \\
& \leq &   {\bf 1}_D F(\cdot, u) \mu.
\end{eqnarray*}
It follows that $(dd^c \varphi)^n = h^n \mu$ where $h^n \in L^1_{\rm loc}(\Omega, \mu)$. By the binomial formula we have
$$
(dd^c (\varphi + v))^n = \sum_{0 \leq k\leq n} C_n^k (dd^c \varphi)^k \wedge (dd^c v)^{n-k}.
$$
By the Monge-Amp\`ere  mixed inequalities we have
$$
(dd^c \varphi)^k \wedge (dd^c v)^{n-k} \geq h^k F(\cdot,v)^{n-k} \mu.
$$
Hence $(dd^c (\varphi + v))^n \geq (h + F(\cdot,v))^n$ in $D$.

On the other hand, since $ {\bf 1}_D (dd^c (\varphi + v))^n \leq {\bf 1}_D (dd^c u)^n \leq F(\cdot,u) \mu$,
we obtain  $h + F(\cdot,v) \leq F(\cdot,u)$ in $D$. Hence  $(dd^c \varphi)^n \leq (F(\cdot,u) - F(\cdot,v))^n \mu$ in $D$.

Since $\varphi \leq 0$ in $\Omega$ we have $u \leq v $ in $D$.
By the hypothesis we have for $z \in D$,
\begin{flalign*}
F(z,u(z)) - F(z,v(z)) &= \int_0^1 \frac{\partial F}{ \partial t}(z, v(z) + t(u(z)-v(z)) (u(z)-v(z)) d s\\
& \leq -\lambda_0 (u(z) - v(z)).
\end{flalign*}
Hence $(dd^c \varphi)^n \leq (-\lambda_0 \varphi)^n \mu$.

If $\varphi \neq 0$, it follows that $E(\varphi) \leq \lambda_0^n \int_\Omega (-\varphi)^{n+} d \mu$.
Hence $\lambda_1 \leq \lambda_0$ by  the minimizing property of  $\lambda_1$ which contradicts our assumption. This prove that $\varphi = 0$ in $\Omega$, hence $u\geq v$ in $\Omega$.

We now prove the existence of solutions by using the fixed point theorem. 
Observe that for $s < 0$ we have
$$
F(z,0) - F(z,s) =  \int_s^0 \frac{ \partial F}{\partial t} (z,t) d t \geq \lambda_0 s.
$$
Applying this inequality to $s=u(z)$, with $u\in \mathcal E^1$, we obtain
$$
F(z,0) - F(z, u(z)) \geq  \lambda_0  u(z),
$$
for $\mu$-a.e. $z \in \Omega$, hence 
\begin{equation}\label{eq: existence F}
F(z,u(z)) \leq F(z,0) + \lambda_0 (- u(z)) \leq f(z) - \lambda_0 u(z)=: g_0 \; \mu\text{-a.e. in }\; \Omega, 
\end{equation}
 Fix a constant $C>0$ and $u\in \mathcal E^1_C$. 
 Then it follows from the inequality \eqref{eq: existence F}  and H\"older inequality, that the measure $\nu := F(u,\cdot)^n \mu$ satisfies the following estimate : there exists a constant $B > 0$ such that for any $\phi \in \mathcal E^1(\Omega)$,
\begin{equation} \label{eq:Ceg0}
\int_\Omega (-\phi) d \nu \leq B \left(\int_{\Omega} (-\phi)^{n+1} d\mu)\right)^{1\slash n+1},
\end{equation}
where $B := (\int_\Omega g_0^{n+1} d \mu)^{n \slash n+1} < + \infty $. On the other hand we have for any $\phi \in \mathcal E^1(\Omega) $,  $\int_{\Omega} (-\phi)^{n+1} d\mu \leq \lambda_1^{-1} E(\phi)$.  
Then for any $\phi \in \mathcal E^1(\Omega)$ we have 
\begin{equation} \label{eq:Ceg}
\int_\Omega (-\phi) d \nu \leq A \left(\int_{\Omega} (-\phi) (dd^c \phi)^n)\right)^{1\slash n+1},
\end{equation}
where $A := B \lambda_0^{-1 \slash n+1}$.
 
 Therefore by \cite[Theorem 8.2]{Ceg98} there exists $v\in \mathcal E^1(\Omega)$ solving the equation 
\[
(dd^c v)^n = F(z,u)^n \mu. 
\]

Using the estimate \eqref{eq: existence F} above, and H\"older's inequality we can write
\begin{flalign*}
	E(v) & = \int_{\Omega} (-v) (dd^c v)^n 
	= \int_{\Omega} F(z,u)^n |v| d\mu  \leq \int_{\Omega}(f-\lambda_0 u)^n |v| d\mu \\
	&\leq \|v\|_{L^{n+1}(\mu)} \|f-\lambda_0 u\|_{L^{n+1}(\mu)}^n\\
	&\leq \|v\|_{L^{n+1}(\mu)} (\|f\|_{L^{n+1}(\mu)}+ \lambda_0 \|u\|_{L^{n+1}(\mu)})^n\\
	&\leq  \lambda_1^{-n/(n+1)} E(v)^{1/(n+1)} \left (A + \lambda_0 \lambda_1^{-n/(n+1)} E(u)^{1/(n+1)}\right)^n,
\end{flalign*}
where $A := \|f\|_{L^{n+1}(\mu)}$.

We have used that $\lambda_1^n\int_{\Omega} (-w)^{n+1}d\mu \leq  E(w)$, for all $w\in \mathcal E^1$.
It then follows that 
\[
E(v)^{1/(n+1)} \leq  A\lambda_1^{-1/(n+1)}  + \lambda_0 \lambda_1^{-1}C^{1/(n+1)}.
\]
Since $\lambda_0 < \lambda_1$, we can choose $C > 1$ so that the right-hand side above is $C^{1/(n+1)}$.  We thus have a map $\eta : \mathcal E^1_C \rightarrow \mathcal E^1_C$ which sends $u$ to $v$ the unique function in $\mathcal E^1$ solving $(dd^c v)^n = F(z,u)\mu$.  Recall that $\mathcal E^1_C$ equipped with the $L^1$-topology is a compact and convex subset of $L^1(\Omega, dV)$. 

We next prove that $\eta$ is continuous. For this, we let $(u_j)$ be a sequence in $\mathcal E^1_C$ converging in $L^1(dV)$ to $u\in \mathcal E^1_C$ and we want to show that $v_j \to v$ in $L^1(dV)$. Since the solution $v$ is uniquely determined by $u$,  it is enough to show that any subsequence of $(v_j)$ converges to $v$.  Extracting a subsequence we can also assume that $u_j$ converges almost everywhere to $u$ with respect to $dV$. Then $z \mapsto g_j(z):= F(z,u_j(z))$ also converges $dV$-almost everywhere to $g(z):=F(z,u(z))$. From \eqref{eq: existence F} and H\"older's inequality, we deduce that 
\[
\|g_j\|_{L^{n+1}(\mu)} \leq \|f + \lambda_0 \vert u_j\vert\|_{L^{n+1}(\mu)} \leq A'. 
\]
We thus infer that $g_j$ converges to $g$ in $L^q(\mu)$, for all $q<n+1$. Extracting once more, we can assume that $g_j \to g$ almost everywhere with respect to $\mu$ and $v_j$ converges almost everywhere and in $L^1(dV)$ to $v$.  Now, define $(\tilde{v}_j=\max_{k\geq j} v_k)^*$. Then $\tilde{v}_j\searrow v$ and 
\[
(dd^c \tilde{v}_j)^n \geq \tilde{g}_j \mu := (\inf_{k\geq j} g_k) \mu.
\]
Letting $j\to +\infty$, we obtain $(dd^c v)^n \geq g \mu$. On the other hand, since $\tilde{v}_j\geq v_j$, we also have 
\[
\int_{\Omega} (dd^c \tilde{v}_j)^n \leq \int_{\Omega} (dd^c v_j)^n = \int_{\Omega} g_j d\mu.  
\]
Letting $j\to +\infty$, we arrive at $\int_{\Omega} (dd^c v)^n \leq \int_{\Omega} gd\mu$, hence $(dd^c v)^n = F(z,u) \mu$. 

The above argument thus shows that $\eta$ is continuous. 
Therefore, by Schauder's fixed point theorem, see  \cite[Theorem B.2, page 302]{Tay11}, there exists a fixed point $u$ of $\eta$, which is a solution to \eqref{eq:DP}.

2. If $\lambda_0 = \lambda_1$ and  $u(z_0) < v(z_0)$ which means that $\varphi \neq 0$ then it follows from  Lemma \ref{lem: supersol is sol}, that $\varphi$ is an eigenfunction associated to the eigenvalue $\lambda_1$.
\end{proof}
As a corollary we obtain Lions type formula.
Consider the following complex Monge-Amp\`ere  equation 
\begin{equation} \label{eq:modifiedMA}
(dd^c u)^n = (1- \lambda u)^n \mu, \, \, u \in \mathcal E^1(\Omega).
\end{equation}

\begin{corollary} \label{cor: Lions formula} Assume $\mu(\Omega)<+\infty$ and  $\mathcal E^1(\Omega) \subset L^{n+1}(\mu)$. For any $\lambda < \lambda_1$ the equation \eqref{eq:modifiedMA} has a unique solution $u_\lambda$. Moreover  we have the following formula:
\begin{equation} \label{eq:LionsFormula}
\lambda_1(\Omega,\mu) = \sup \{\lambda > 0 ; \eqref{eq:modifiedMA} \, \, \text{has a subsolution}\}\cdot
\end{equation}

\end{corollary}
\begin{proof} Let $\gamma_1$ denotes the right hand side of  the identity \eqref{eq:LionsFormula}. By Theorem \ref{thm:Existence-Uniqueness} for any $0 < \lambda < \lambda_1 := \lambda_1(\Omega,\mu)$ the equation \eqref{eq:modifiedMA} has a (unique) solution
$u \in \mathcal E^1(\Omega)$. 
Hence $ \gamma_1 \geq \lambda_1$. On the other hand, if $\lambda>0$ is such that \eqref{eq:modifiedMA} has a subsolution $u\in \mathcal E^1(\Omega)\setminus\{0\}$, then
\[
(dd^c u)^n \geq (1-\lambda u)^n \mu \geq (-\lambda u)^n \mu.
\]  
Thus, by  Theorem \ref{thm: uniqueness of the eigenvalue},  we have $\lambda \leq \lambda_1$, and the result follows. 
\end{proof}

 A well-known fact concerning the first eigenvalue of a linear elliptic operator is that all corresponding eigenfunctions have constant sign. In particular, any negative eigenfunction of the Laplace operator is subharmonic. In the following result, using the envelope technique, we prove a similar phenomenon for the Monge-Amp\`ere operator: any smooth negative eigenfunction is automatically  plurisubharmonic. 
\begin{prop}
	Assume $\mu= fdV$, $f>0$ is smooth, and $v \in C^2(\Omega) \cap C^0(\bar \Omega)$ is a negative function which vanishes on $\partial \Omega$, and it satisfies
	\[
	(dd^c v)^n \leq  (-\lambda_1 v)^{n} fdV,
	\]
	pointwise in $\Omega$. Then $v$ is plurisubharmonic in $\Omega$.  
\end{prop}
Note that, in the above result we do not assume $dd^c v\geq 0$. 
\begin{proof}
	Let $u=P(v)\in \mathcal E^1(\Omega)$. Then $u$ is a supersolution. Invoking Lemma \ref{lem: supersol is sol} we see that $u$ is a solution. The measure $(dd^c u)^n =(-\lambda_1 u)^n fdV$ does not charge the open set $\{u<v\}$. Since $u<0$ and $f>0$, it follows that $\{u<v\}$ has zero Lebesgue measure hence it is empty.  This also implies that $v$ is plurisubharmonic. 
\end{proof}
 \subsection{More general boundary values}
 In this section we assume that $\mu(\Omega)<+\infty$ and $\mathcal E^1(\Omega) \subset L^{n+1}(\Omega,\mu)$.

Let  $h \in  C^0(\partial \Omega)$ a boundary data. Assume that $h$ extends as a plurisubharmonic function $H \in \PSH(\Omega) \cap  C^0(\partial \Omega) $.
Let $\Phi $ be  the Perron-Bremermann envelope 
of the family of plurisubharmonic functions $u \in \PSH(\Omega)$ such that $u^* \leq h$ in $\partial \Omega$.
Then by  \cite{BT76} $\Phi \in \PSH(\Omega) \cap C^0(\bar \Omega)$ is the solution to the homogeneous complex Monge-Amp\`ere  equation :
 $$
 (dd^c \Phi)^n = 0, \, \, \text{on} \, \, \Omega, \, \, \Phi_{\mid \partial \Omega} = h.
 $$

 Following Cegrell, we denote by $\mathcal E^1(\Omega,\Phi)$ the set of $u \in \PSH(\Omega)$ such that there exists $v \in \mathcal E^1(\Omega)$ with $v + \Phi \leq u \leq \Phi$. In other words this means that
 $\varphi := P_\Omega (u-\Phi,0) \in \mathcal E^1(\Omega)$.
 
 We now consider the following complex Monge-Amp\`ere  equation with boundary values $\Phi \in \mathcal {MPSH}(\Omega)$:
\begin{equation} \label{eq:DP2}
	(dd^c u)^n  -  F(\cdot,u)^n \mu = 0, \;  \; u \in \mathcal E^1(\Omega,\Phi), 
\end{equation}
where $F :  \Omega \times \R \longrightarrow \R^+$ is a Borel function such that 
\begin{itemize}
\item the density $z \longmapsto f(z) := F(z,\Phi(z)) $ satisfies  $f \in L^{n+1}(\Omega,\mu)$, and 
\item for $\mu$-a.e. $z \in \Omega$, $t \longmapsto F(z,t)$ is continuous in $]-\infty, \max_{\partial \Omega} \Phi]$. 
\end{itemize}

 Then we have the following general result.
\begin{theorem}  \label{thm:Existence-Uniqueness2} Assume that there exists $ \lambda_0 > 0 $ such that
$$
\frac{\partial F}{\partial u} \geq - \lambda_0 , \, \, \mathrm{in} \, \, \Omega \times  \R.
$$ 
Then  

(1) if $\lambda_0 < \lambda_1$,  the maximum principle holds for the operator 
$(dd^c u)^n - F(\cdot,u) \, \mu$ on $\mathcal E^1(\Omega,\Phi)$  i.e. if $u \in \mathcal E^1(\Omega,\Phi)$ is a supersolution and $v \in \mathcal E^1(\Omega,\Phi)$ is a subsolution then $u \geq v$ in $\Omega$. 

Moreover the Dirichlet problem \eqref{eq:DP2} admits a unique solution $\varphi \in \mathcal E^1(\Omega,\Phi)$. In particular, if $F(z,\Phi(z))\equiv 0$, the unique solution is  $\Phi$.

(2) If $\lambda_0 = \lambda_1$ and $u, v \in \mathcal E^1(\Omega,\Phi)$ are two solutions to \eqref{eq:DP2}  such that $u(z_0) < v(z_0)$ for some point $z_0 \in \Omega$, then $\varphi := P(u-v,0)<0$ is an eigenfunction associated  to the eigenvalue $\lambda_1$.

\end{theorem}
\begin{proof} The proof of the comparison principle goes in the same way as in the previous theorem.
This implies uniqueness of the solution.

To prove the existence we first observe that by the same reasoning which leaded to \eqref{eq: existence F}, we see that for any $v \in \mathcal E^1(\Omega,\Phi)$ we have for $\mu$-almost every $z \in \Omega$, 
 \begin{equation}\label{eq: inequality F}
F(z,v(z)) \leq F(z,\Phi(z)) + \lambda_0 (\Phi(z) - v(z)) \leq g (z)  \; \text{ in }\; \Omega. 
\end{equation}
where  $g := f + \lambda_0( \Phi-v) \in L^{n +1} (\Omega,\mu)$.

 By Theorem \ref{thm:Existence-Uniqueness}, there exists $u_0 \in  \mathcal E^1(\Omega)$ such that
$$
(dd^c u_0)^n = (f- \lambda_0 u_0)^n \mu.
$$
Then $v_\Phi := u_0 + \Phi \in  \mathcal E^1(\Omega,\Phi)$ is a subsolution to the Dirichlet problem \eqref{eq:DP2}. Indeed  by \eqref{eq: inequality F} we have
$$
(dd^c v_\Phi)^n \geq (dd^c u_0)^n = (f + \lambda_0(\Phi - v_\Phi))^n \mu \geq F(\cdot,v_\Phi)^n \mu.
$$

We will  use the Schauder fixed point theorem again. Consider the  following convex set 
$$
\mathcal C := \{ v \in PSH (\Omega) ; v_\Phi \leq v \leq \Phi\}\cdot
$$
  
 Then $\mathcal C \subset \mathcal E^1(\Omega,\Phi)$  is a non empty convex and compact set for the $L^1(\Omega)$-topology. 
 
 Let $v \in \mathcal C$.
 Observe that $ 0 \leq  f + \lambda_0(\Phi- v) \leq f - \lambda_0 u_0 =: g_0 \in L^{n +1} (\Omega,\mu)$.

  Then it follows from the inequality \eqref{eq: inequality F}  and H\"older inequality, that the measure $\nu := F(v,\cdot)^n \mu$ satisfies the estimate \eqref{eq:Ceg}.
 
Therefore by \cite[Theorem 8.2]{Ceg98} there exists a unique solution $u \in  \mathcal E^1(\Omega,\Phi)$ to the equation 
 \begin{equation} \label{eq:DP3}
 (dd^c u)^n = F(v,\cdot)^n  \mu,
 \end{equation} 
 on $\Omega$.

Since $v \geq v_\Phi$, by  the inequality \eqref{eq: inequality F} we have
$$
F(\cdot,v) \leq f + \lambda_0(\Phi - v) \leq f + \lambda_0(\Phi - v_\Phi) = f - \lambda_0 u_0.
$$
Hence $(dd^c u)^n \leq (dd^c v_\Phi)^n$ on $\Omega$. By the comparison principle it implies that $v_\Phi \leq u \leq \Phi$, hence $u  \in \mathcal C$.

It remains to show that $T : \mathcal C \longrightarrow \mathcal C$ is continuous. Indeed let $(v_j)$ be a sequence in $\mathcal C$ which converges to $v \in \mathcal C$ in the $L^1(\Omega)$-topology.
Then we know that $v = (\limsup_{j \to +\infty} v_j)^*$ in $\Omega$. By definition for any $j \in \N$, $u_j := T(v_j) \in  \mathcal E^1(\Omega,\Phi)$ solves the following equation
$$
(dd^c u_j)^n = F(\cdot,v_j) \mu.
$$
We want to show that $u_j \to u$, where $u := T(v)$ satisfies  $(dd^c u)^n = F(\cdot,v) \mu$.
Up to taking a subsequence we can assume that $u_j \to U \in \mathcal C$ in $L^1(\Omega)$ and almost everywhere in $\Omega$ for the Lebesgue measure. 

By uniqueness it's enough to prove that $(dd^c u_j)^n \to (dd^c U)^n$  weakly on $\Omega$. Indeed this will imply that $U = u = T(v)$, proving the continuity of $T$.

To prove the weak convergence $(dd^c u_j)^n \to (dd^c U)^n$, we first observe that  $(dd^c u_j)^n \leq (f - \lambda_0)^n) \mu$ for any $j$.
Then 
$$
\int_\Omega \vert u_j -U\vert (dd^c u_j)^n \leq \int_\Omega \vert u_j -U\vert  (f - \lambda_0 u_0)^n d \mu.
$$
Since $\mu$ vanishes on pluripolar sets, it follows from \cite[Corollary 1.5]{CK06} that $u_j \to U$ $\mu$-a.e. in $\Omega$.
Moreover since 
$$
\vert u_j -U\vert  (f - \lambda_0 u_0)^n \leq (-u_0)((f - \lambda_0 u_0)^n \in L^{n+1}(\Omega,\mu),
$$ 
it follows from Lebesgue convergence theorem that 
$\lim_{j \to + \infty} \int_\Omega \vert u_j -U\vert (dd^c u_j)^n =0$.
A classical argument  as in \cite[Lemma 5.2]{Ceg98}  yields the convergence $(dd^c u_j)^n \to (dd^c U)^n$ in the sense of Radon measures on $\Omega$.
 \end{proof}

\section{The iterative approach}\label{sect: iterative}
We assume in this section that $\mu$ is a positive Borel measure on $\Omega$ such that $\mu(\Omega) >0$ and $\mathcal E^1(\Omega) \subset L^{n+1} (\Omega)$ so that $I_{\mu}$ is finite on $\mathcal E^1(\Omega)$.

The goal of this section is to prove Theorem \ref{thm: iterative intro}.

We already observed that $I_\mu (u) > 0$ if $ u \in \mathcal E^1(\Omega)\setminus \{0\}$. 
 
 We fix $u_0\in \mathcal E^1(\Omega)\setminus \{0\}$. We define inductively a sequence $(u_k)_{k \in \N}$ in $\mathcal E^1(\Omega)\setminus \{0\}$ as follows :  for each $k \in \N$,
\[
(dd^c u_{k+1})^n = R(u_k) (-u_k)^n \mu, \, \, \text{where} \, \, \; R(u_k):= \frac{E(u_k)}{I_{\mu}(u_k)}. 
\]
Since $I_{\mu}$ is finite on $\mathcal E^1(\Omega)$, an application of H\"older's inequality yields, for all $v\in \mathcal E^1(\Omega)$, 
\[
\int_{\Omega}(-v) (-u_k)^n d\mu \leq \left ( \int_{\Omega}(-v)^{n+1} d\mu  \right)^{1/(n+1)} \left ( \int_{\Omega}(-u_k)^{n+1} d\mu  \right)^{n/(n+1)}<+\infty. 
\] 
The existence of $u_{k+1}$ thus follows from \cite{Ceg98}. 

It follows from \cite[Lemma 3.4]{Zer25} that 
\begin{equation}
	\label{eq: energy monotone}
	\frac{E(u_{k+1})}{I_{\mu}(u_{k+1})^{1/(n+1)}}  \leq \frac{E(u_k)}{I_{\mu}(u_k)^{1/(n+1)}},
\end{equation}
and 
\begin{equation}
	\label{eq: energy monotone 2}
R(u_{k+1})I_{\mu}(u_{k+1})^{n/(n+1)} \leq R(u_k) I_{\mu}(u_k)^{n/(n+1)}. 
\end{equation}
From \eqref{eq: energy monotone} we thus have  
\[
\frac{E(u_{k})}{I_{\mu}(u_{k})^{1/(n+1)}}  \leq \frac{E(u_0)}{I_{\mu}(u_0)^{1/(n+1)}} =C_0. 
\]
On the other hand, we also have $\lambda_1(\mu) I_{\mu}(u_k) \leq E(u_k)$. Therefore, 
\[
\lambda_1(\mu)^{1/(n+1)}E(u_k) \leq C_0E(u_k)^{1/(n+1)},
\]
and we infer $E(u_k) \leq \lambda_1(\mu)^{1/n} C_0^{1+1/n}$ is uniformly bounded. From this we also infer that $I_{\mu}(u_k)$ is uniformly bounded from above.

It is quite surprising that the following quantities are monotone along the iterative scheme. 
\begin{lemma}\label{lem: monotone}
	Along the above iterative scheme the sequences $(E(u_k))_{k \in \N}$ and $(I_{\mu}(u_k))_{k \in \N}$ are both increasing, while $(R(u_k))_{k \in \N}$ is decreasing in $k$. 
\end{lemma}

\begin{proof}
We provide two methods to prove that the sequence $(E(u_k))_{k \in \N}$ is increasing. The first approach using Cegrell's inequality is inspired by Q.N. L\^e,  \cite{LeQN20}, who studied the eigenvalue problem for real Monge-Amp\`ere equations. Although the analogous Cegrell's inequality goes in the opposite direction, the main idea can be recycled in our seting. 

Using Cegrell's inequality \eqref{eq: Cegrell inequality 3} for $u=u_k$ and $v=u_{k+1}$, we obtain 
\begin{flalign*}
	\int_{\Omega} R(u_k) (-u_k)^{n+1} d\mu &= \int_{\Omega} (-u_k) (dd^c u_{k+1})^n  d\mu \leq  E(u_k)^{1/(n+1)} E(u_{k+1})^{n/(n+1)}. 
\end{flalign*}
Recalling the definition of $R(u_k)$, we thus get 
\[
E(u_k) \leq  E(u_k)^{1/(n+1)} E(u_{k+1})^{n/(n+1)},
\]
which yields $E(u_k)\leq E(u_{k+1})$. 

Alternatively, one can prove that $E(u_k)$ is increasing in $k$ by using the convexity of $E$ along affine curves in $\mathcal E^1(\Omega)$. Indeed, from the definition of $u_k$ we have 
\begin{eqnarray} \label{eq:energy relation}
\int_{\Omega} (u_k - u_{k+1}) (dd^c u_{k+1})^n & =  E(u_{k+1}) - R(u_k) I_{\mu}(u_k)  \nonumber\\
&= E(u_{k+1})-E(u_k).
\end{eqnarray}
On the other hand, we have 
\begin{eqnarray}\label{eq: energy convexity}
 E(u_{k+1})-E(u_k) &= & \sum_{j=0}^{n} \int_{\Omega} (u_k-u_{k+1}) (dd^c u_k)^j \wedge (dd^c u_{k+1})^{n-j}
  \nonumber   \\
&=& \int_{\Omega} (u_k-u_{k+1}) (dd^c u_{k+1})^n  + A,
\end{eqnarray}
where $A := \sum_{j = 1}^{n} \int_{\Omega} (u_k-u_{k+1}) (dd^c u_k)^j \wedge (dd^c u_{k+1})^{n-j}.$
 From  \eqref{eq:energy relation} and \eqref{eq: energy convexity} 
we infer that $A=0$.

By the formula \eqref{eq:Integration-part} we have for $1 \leq j \leq n$
$$
\int_{\Omega} (u_k-u_{k+1}) (dd^c u_k)^j \wedge (dd^c u_{k+1})^{n-j} \leq \int_{\Omega} (u_k-u_{k+1}) (dd^c u_{k+1})^n.
$$ 
Hence 
 $A \leq n\int_{\Omega} (u_k-u_{k+1}) (dd^c u_{k+1})^n $.  It thus follows that $E(u_{k+1})-E(u_k) \geq 0$.

Since $E(u_k)$ is increasing in $k$, using \eqref{eq: energy monotone} we deduce  that $I_{\mu}(u_k)$ is also increasing in $k$. From \eqref{eq: energy monotone 2} we then infer that $R(u_k)$ is decreasing in $k$. 
\end{proof}

\begin{theorem}\label{thm: existence iterative}
		If $I_{\mu}$ is continuous on each $\mathcal E^1_C$ then the sequence $(u_k)$ converges in $L^1$ to an eigenfunction $u$.  
\end{theorem}
\begin{proof}
		Let $(u_{k(j)})$ be a subsequence of $(u_k)$, converging in $L^1(dV)$ and almost everywhere to some $u\in \mathcal E^1(\Omega)$. 
It follows from Lemma \ref{lem: convergence} (and the continuity assumption on $I_{\mu}$) that 
	\[
	\int_{\Omega} (u_{k(j)}-u) (dd^c u_{k(j)})^n = \int_{\Omega} (u_{k(j)}-u) (-u_{k(j)-1})^{n} d\mu \to 0. 
	\]
	By Fatou's lemma we have 
	\[
	\limsup_j \int_{\Omega} (u_{k(j)}-u) (dd^c u)^n \leq 0.
	\]
	By the inequality \eqref{eq:Integration-part} we also have 
	\[
	\int_{\Omega} (u_{k(j)}-u) (dd^c u)^n \geq \int_{\Omega} (u_{k(j)}-u) (dd^c u_{k(j)})^n. 
	\]
	It thus follows that $\int_{\Omega} (u_{k(j)}-u) (dd^c u)^n$ also converges to $0$. Therefore, the inequality \eqref{eq:Integration-part} and the cocycle formula yield $E(u_j) \to E(u)$. 
	
	Since  $I_{\mu}(u_{k(j)}) \nearrow I_{\mu}(u)$, we also have $R(u_{k(j)}) \searrow R(u)$. Setting $v=\lim_{j} u_{k(j)+1}$ and arguing as above we also get 
	\[
	E(u_{k(j)}) \nearrow E(v)=E(u), \; I_{\mu}(u_{k(j)}) \nearrow I_{\mu}(u), \; R(u_{k(j)}) \searrow R(v)=R(u). 
	\]
	Setting $\tilde{v}_j:= (\sup_{l\geq j} u_{k(l)})^*$, and using the pluripotential maximum principle, we arrive at 
	\[
	(dd^c \tilde{v}_j)^n \geq R(v) \left(\inf_{l\geq j} |u_{k(l)}|\right)^n \mu. 
	\]
	Letting $j\to +\infty$, we obtain 
	\[
	(dd^c v)^n \geq  R(u) (-u)^n \mu. 
	\]
	Using Cegrell's inequality as in the proof of Lemma \ref{lem: monotone} we have 
	\begin{flalign*}
		E(u) & = R(u) \int_{\Omega} (-u)^{n+1} d\mu \leq \int_{\Omega} (-u) (dd^c v)^n \\
		& \leq E(u)^{1/(n+1)} E(v)^{n/(n+1)} = E(u).  
	\end{flalign*}
	Thus the two measures $(-u)(dd^c v)^n \geq  R(u) (-u)^{n+1} \mu$ have the same total mass over $\Omega$. They must be equal, hence (since $u<0$), we have the equality $(dd^c v)^n = R(u) (-u)^n \mu$. 
	Having this we can proceed as in \cite{Zer25} to conclude that $v=cu$ for some $c>0$. Then $c=1$ because $E(u)=E(v)$. To prove that any other cluster point of $(u_k)$ coincides with $u$, we repeat the arguments in \cite{Zer25} and we use the uniqueness result, Corollary \ref{cor: unique}.
\end{proof}

We end this section with a proof of Theorem \ref{thm: uniqueness of the eigenvalue} using the iterative method.

\begin{proof}[Proof of Theorem \ref{thm: uniqueness of the eigenvalue} via the iterative method]
	Assume $(dd^c \psi)^n \geq (-\lambda \psi)^n \mu$ with $\psi \in \mathcal E^1(\Omega) \setminus \{0\}$. This implies in particular that  $I_{\mu}(\psi)<+\infty$. Our goal is to prove that $\lambda\leq \lambda_1$. Recall that $\lambda_1$ is defined as 
	\[
	\lambda_1:= \inf \left \{ \frac{E(u)}{I_{\mu}(u)} \; : \; u\in \mathcal E^1, \; \psi \leq u<0\right \},
	\] 
	and it does not depend on $\psi$, as shown in Lemma \ref{lem: def of lambda 1}. Let us emphasize here that we do not assume $I_{\mu}$ is finite on $\mathcal E^1$ but merely that $I_{\mu}(\psi) <+\infty$. 
	
	Assume by contradiction that $\lambda>\lambda_1$, and take $u\in \mathcal E^1(\Omega)$, $u\neq 0$, such that $R(u)<\lambda^n$. For $t$ sufficiently large we have $R(\max(u,t\psi))< \lambda$. Replacing $\psi$ by $t\psi$,  we can assume $t=1$. Consider the iterative scheme starting at $u_0:=\max(u,\psi)$. We prove that there is a subsequence converging to an eigenfunction of \eqref{eq: MA eigenvalue}. Recall that $u_{k}$ is defined inductively as the unique solution in $\mathcal E^1$ to the Monge-Amp\`ere equation $(dd^c u_k)^n =R(u_{k-1}) (-u_{k-1})^n \mu$.  Since 
\[
R(u_0) (-u_0)^n \leq (-\lambda \psi)^n\leq (dd^c \psi)^n,
\] 
it follows from the comparison principle that $u_1 \geq \psi$. Since $R(u_1)\leq R(u_0) \leq \lambda$ we have
\[
R(u_1) (-u_1)^n \leq (-\lambda \psi)^n\leq (dd^c \psi)^n.
\] 
The comparison principle again gives $u_2\geq \psi$. 
  Since $R(u_k)$ is decreasing in $k$, we can prove by induction that $u_k \geq \psi$ for all $k$.  Since the sequence $u_k$ stays in a compact set of $L^1(dV)$, we can find a subsequence $v_j:= u_{k(j)}$ converging to $\varphi\in \mathcal E^1$ in $L^1(dV)$. Extracting once more, with the new sequence still denoted by $v_j$, we can assume that the convergence holds in  $L^1(dV)$, $L^1(\mu)$, almost everywhere with respect to $dV$ as well as $d\mu$.    We can thus use the dominated convergence theorem to show that $E(v_j)$, $I_{\mu}(v_j)$ converge to $E(\varphi)$, $I_{\mu}(\varphi)$. Arguing as in the proof of Theorem \ref{thm: existence iterative}, we can show that $\varphi$ is an eigenfunction.  
  
Now, we have 
\[
(dd^c \varphi)^n = R(\varphi) (-\varphi)^n \mu, \; (dd^c \psi)^n \geq \lambda^n (-\psi)^n \mu.
\]
Thus, for $a= R(\varphi)^{-1} \lambda^n >1$, we have $(dd^c a \varphi)^n \leq (dd^c \psi)^n$ on $\{a\varphi < \psi\}$.
By the domination principle, we have $a\varphi \geq \psi$. Repeating the same argument we get $a^m \varphi \geq \psi$, for all $m>1$, which is impossible. 
\end{proof}

\section{The eigenvalue problem for complex Hessian operators} \label{sect: Hessian}
Given $1\leq m\leq n$, the complex Hessian operator $H_m(u)$, acting on $m$-subharmonic functions, is defined as 
\[
H_m(u)= (dd^c u)^m \wedge \beta^{n-m}.
\]
We refer the reader to \cite{Bl05} for background on $m$-subharmonic functions and to \cite{BZ24} for the variational approach to the eigenvalue problem associated with complex $m$-Hessian operator. We are looking for pairs $(\lambda,u)$ solving
\begin{equation}
	\label{eq: eigen Hes}
	H_m(u) = (-\lambda u)^m \mu,
\end{equation} 
where $\lambda>0$ and $u$ is an $m$-subharmonic function in the Cegrell's finite energy class $\mathcal E^1_m(\Omega)\setminus \{0\}$, and $\mu$ is a positive Borel measure vanishing on $m$-polar sets. Here $\Omega$ is a bounded $m$-hyperconvex domain in $\mathbb C^n$. 

Our method using envelopes can be extended to the Hessian setting, giving analogous results. In particular, open questions raised at the end of \cite{BZ24} have positive answers. 
\begin{theorem}\label{thm: eigen Hess unique}{\; }
\begin{enumerate}
	\item If  $(\lambda,\varphi)$ is a subsolution to \eqref{eq: eigen Hes} then $\lambda\leq \lambda_{1,m}$.  
	\item If $(\lambda,\varphi)$ solves \eqref{eq: eigen Hes} then $\lambda = \lambda_{1,m}$. 
	\item If $(\lambda, \varphi)$ solves \eqref{eq: eigen Hes} while $(\lambda,\psi)$ is a subsolution then $(\lambda,\psi)$ is actually a solution and $\psi=c\varphi$ for some constant $c>0$. 
	\item If $I_{\mu}$ is continuous on each $\mathcal E^1_{m,C}(\Omega)$ then the eigenvalue problem has a unique solution $(\lambda_{1,m},u)$. 
 Moreover, the iterative scheme, defined by 
	\[
	H_m(u_{k+1}) = R(u_k) (-u_k)^m \mu, \; u_0 \in \mathcal E^1_m(\Omega)\setminus \{0\}, 
	\]
	produces a sequence $(\lambda_j, u_j)$ converging to $(\lambda_1,u)$, where $u$ is an eigenfunction. 
\end{enumerate}
\end{theorem}
Here $\lambda_{1,m}$ is defined as 
\[
\lambda_{1,m} = \inf \left \{ \frac{\int_{\Omega} (-u) H_m(u)}{\int_{\Omega} (-u)^{m+1} d\mu} \; : \; \varphi \leq u < 0, \; u \in SH_m(\Omega) \right \}.
\]
The existence of a subsolution $\varphi$ implies $\int_{\Omega} (-\varphi)^{m+1} d\mu <+\infty$, hence $\lambda_{1,m}$ is finite. Arguing as in the Monge-Amp\`ere case, we can show that $\lambda_{1,m}$ does not depend on $\varphi$. 
\begin{proof}
The fine properties of $m$-subharmonic envelopes required in the preceding sections were established in \cite{ACLR24}. The entire proof can therefore be adapted to the complex Hessian setting, and the details are left to the reader.
\end{proof}

The analogous version of Theorem \ref{thm:Existence-Uniqueness2} also holds in the Hessian setting. We leave the detail to the reader. We end this section with the following regularity result for the eigenfunctions. 

\begin{theorem}\label{thm: eigen Hessian bounded}
	If $\mu= fdV$ for some $f\in L^p(dV)$ with $p>n/m$, then the eigenfunction $u$ is  continuous on $\overline{\Omega}$. If $\Omega$ is strongly $m$-pseudoconvex then $u$ is H\"older continuous on $\overline{\Omega}$.
\end{theorem}
\begin{proof}
	The existence of the eigenvalue and eigenfunctions follows from \cite{BZ24}. 
	To prove that $u$ is bounded, the argument in \cite{BZ24} is to show that $|u|^mf$ is in $L^q$ for some $q>n/m$, and this is possible under some condition on $m$ and $p$ (when $m=n$ this holds for all $p>1$ since psh functions have strong integrability properties). We will remove the conditions in \cite[Theorem 1.1]{BZ24} by combining the method of Lions \cite{Lio85} and our uniqueness result. 
	
	Define $\lambda_0=\sup S$, where $S$ is the set of all $\lambda>0$ such that there is a solution $u\in \mathcal E^1_m$ to 
	\begin{equation}\label{eq: Lions Hess}
		H_m(u) = (1-\lambda u)^m fdV.
	\end{equation}  
	We first prove that $S$ is not empty. By \cite{DK14}, we can find a bounded $u_0\in SH_m(\Omega)$ solving $(dd^c u_0)^n =fdV$. Then the function $
	u_1:= Cu_0$ satisfies 
	\[
	H_m(u_1) = C^m fdV \geq (1-\lambda C u_0)^m fdV = (1-\lambda u_1)^m fdV,
	\]
	if $C \geq 1+\lambda C  \|u_0\|_{\infty}$. Thus for small $\lambda$, there is a subsolution to  \eqref{eq: Lions Hess}. We can then inductively define $u_k \in \mathcal E^1_m(\Omega)$, $k\geq 2$, to be the unique $m$-sh solution to \eqref{eq: Lions Hess} with $u_{k-1}$ on the right-hand side. The comparison principle yields that the sequence $(u_k)$ is increasing in $k$, hence it converges to a solution to \eqref{eq: Lions Hess} which cannot be zero. 
	
	Thus $S$ is not empty. Observe that if $\lambda\in S$ then, for some $u\in \mathcal E^1_m\setminus \{0\}$, $(\lambda, u)$ is a subsolution to \eqref{eq: eigen Hes}. Therefore $\lambda \leq \lambda_{1,m}$, hence $\lambda_0\leq \lambda_1$. Arguing as in Corollary \ref{cor: Lions formula} we see that $\lambda_0=\lambda_1$. Let $(\lambda_j)$, $(u_j)$ be such that 
	\[
	H_m(u_j) = (1-\lambda_j u_j)^m fdV, \; \lambda_j\nearrow \lambda_1. 
	\]
	If $\sup_{\Omega} |u_j|$ is bounded then the right-hand side above is uniformly in $L^p$, $p>n/m$, thus by \cite[Theorem 2.5]{DK14}, $u_j$ uniformly converges to a bounded solution to $H_m(u)= (1-\lambda_1 u)^m fdV$. Then $(\lambda_1,u)$ is a subsolution to \eqref{eq: eigen Hes}, and by Theorem \ref{thm: eigen Hess unique}, $u$ must be an eigenfunction, $(1-\lambda_1 u)^m fdV = (-\lambda_1 u)^m fdV$, which is impossible. Therefore, $\sup_{\Omega}|u_j| \to +\infty$. Now, for $v_j:= \|u_j\|_{\infty}^{-1}u_j$, we have 
	\[
	H_m(v_j) = \left( \|u_j\|_{\infty}^{-1} - \lambda_j v_j \right)^m fdV. 
	\]
	Since the right-hand side is uniformly in $L^p$, $p>n/m$, we infer by \cite[Theorem 2.5]{DK14} that $v_j$ uniformly converges to $v\in SH_m(\Omega)$, and $v$ solves \eqref{eq: eigen Hes}. By construction $v$ is bounded since $\|v_j\|_{\infty}=1$, hence by \cite[Theorem 2.5]{DK14} $u$ is continuous on $\overline{\Omega}$. Our uniqueness theorem then shows that the eigenfunctions of \eqref{eq: eigen Hes} are bounded. They are even H\"older continuous when $\Omega$ is strongly $m$-pseudoconvex, see \cite{BZ20}, \cite{CZ24}, and \cite{KN25}. 
\end{proof}

\section{The eigenvalue problem for the real Monge-Amp\`ere  operator}
It is possible to use the envelope method to study the eigenvalue problem for the real Monge-Amp\`ere operator. However to avoid repetition we will show that the eigenvalue problem for the real Monge-Amp\`ere  operator can be reduced to the eigenvalue problem for the complex Monge-Amp\`ere operator.  

\subsection{The real Monge-Amp\`ere  operator}
Let $D \Subset \R^n$ be a bounded {\it convex} domain. We denote by ${\rm CV}(D)$ the positive cone of real valued convex functions in $D$, and by ${\rm CV}_0(D)$ the positive cone of  real valued convex functions in $D$ that are continuous on $\bar D$ and vanish on the boundary $\partial D$.  

The following result is well-known but we include a proof for completeness. 
\begin{lemma}\label{lem: distance concave}
	The function $x \mapsto d(x):= {\rm dist}(x,\mathbb R^n\setminus D)$ is concave in $D$, and it vanishes on $\partial D$. 
\end{lemma}
\begin{proof}
	Let $x_0,x_1\in D$ and $x_t=(1-t)x_0+tx_1\in D$. We want to prove that $d(x_t) \geq  r_t:=(1-t)d(x_0)+t d(x_1)$. Fix $y\in \mathbb R^n$ such that $\|y-x_t\|<r_t$. It suffices to prove that $y\in D$. Observe that we can write 
	\[
	y = (1-t) y_0 +t y_1, \; y_j:= x_j + (y-x_t) \frac{d(x_j)}{r_t}, \; j=0,1. 
	\] 
	Now, using $\|y_j-x_j\| \leq d(x_j)$, we arrive at $y_j\in D$, $j=0,1$, hence $y\in D$ by convexity of $D$. 
\end{proof}

If $v \in {\rm CV}(D) \cap C^2(D)$ we define the real Monge-Amp\`ere  measure of $u$ as 
$$
\text{M}_\R(u) := \text{det} \, \left(\frac{\partial^2 u}{\partial x_j \partial x_k}\right) dV(x),
$$
where $dV(x) = \bigwedge_{1 \leq j\leq n} dx_j$ is the euclidean volume form on $D$ identified with the $n$-dimensional Lebesgue measure on $D$.

An important result of Alexandrov \cite{AZ67} shows that  this definition can be  extended to any convex function $u \in {\rm CV}(D)$ so that $\text{M}_\R(u)$ is a Borel measure on $D$ defined as follows: 
 its mass on any Borel set $K\subset D$ is the euclidean volume of the gradient image of $K$ i.e.
$$
\int_K \text{M}_\R(u) = \vert \nabla u (K)\vert_n,
$$
where $ \vert \cdot \vert_n$ is the $n$-dimensional volume and $\nabla u (K) := \bigcup_{a\in K} \nabla u (a),$ where
$$
  \nabla u (a) := \{ \xi \in \R^n ; u(x) \geq u(a) + \langle \xi, x-a\rangle, \forall x \in D\}. 
$$ 
There is a deep connection between the real Monge-Amp\`ere  operator and the complex Monge-Amp\`ere  operator. To describe it we introduce the logarithmic map
$L : (\C^*)^n \longrightarrow \R^n$ defined for $z = (z_1, \cdots, z_n) \in  (\C^*)^n$ by $x = L (z)$, where 
$$
 L (z) := (\log\vert z_1\vert, \cdot, \log \vert z_n\vert).
$$
It is a smooth and proper map from $(\C^*)^n$ onto $\R^n$, so the set 
$$
\Omega := L^{-1}(D) \Subset (\C^*)^n,
$$
is a bounded hyperconvex domain invariant by the action of the torus $\mathbb T^n := (\mathbb S^1)^n$ such that $ L(\Omega) = D$. Indeed, to see that $\Omega$ is hyperconvex, we can consider the function $\rho(z) = - {\rm dist}(L(z), \mathbb R^n\setminus D)$. By Lemma \ref{lem: distance concave} above the function $x\mapsto -{\rm dist}(x,\mathbb R^n \setminus D)$ is convex in $D$ hence $\rho$ is a plurisubharmonic exhaustion function of $\Omega$. 

For each $\zeta :=(\zeta_1, \cdots, \zeta_n) \in \mathbb T^n $,  the map 
$$
z \in (\C^*)^n \longmapsto \tau_\zeta (z) := (\zeta_1 z_1, \cdots, \zeta_n z_n)  \in (\C^*)^n
$$ 
is a biholomorphism such that $\tau_\zeta(\Omega) = \Omega$ and $\vert J \tau_\zeta (z)\vert \equiv 1$ in $(\C^*)^n$ .

Recall that a function $\psi : \Omega \longrightarrow \R$ is toric if $\psi(z) = \psi (\vert z_1\vert, \cdot,\vert z_n\vert)$ for any $z \in \Omega$ which means that it's invariant by the action of the torus $\mathbb T^n$ on $\Omega$. This means that for any $\zeta \in \mathbb T^n$, $ \psi \circ \tau_\zeta = \psi$ in $\Omega$.

It is well known that if a function $u : D \longrightarrow \R$  is convex in $D$ then $\varphi = \varphi_L := u \circ L$ is a (continuous) toric plurisubharmonic in $\Omega$. Conversely if $\varphi$ is a toric plurisubharmonic function in $\Omega$ then the function $u(x) := \varphi (e^{x_1}, \cdots, e^{x_n})$ is a convex function in $D$ such that $\varphi = u \circ L$ in $\Omega$.

Moreover the complex Monge-Amp\`ere  measure of $\varphi = u \circ L$ is related to the real Monge-Amp\`ere  measure of $u$ by the following formula (see \cite{BB13} and \cite{CGSZ19}) : for any continuous test function $\chi$ in $D$, 
\begin{equation} \label{eq:MAR-MAC}
 \int_\Omega \chi \circ L (dd^c \varphi)^n = \int_D \chi \mathrm{M}_\R(u).
\end{equation}
This means that $L_*((dd^c \varphi)^n) = M_\R(u)$ in the sense of Borel measures on $D$. 
Observe that $\mu := (dd^c \varphi)^n$ is a toric measure on $\Omega$ which means that it is invariant by the torus action on $\Omega$ i.e. for any $\zeta \in \mathbb T^n$, 
$(\tau_\zeta)_* \mu = \mu$ in the sense of Borel measures on $\Omega$.

\smallskip

We will need the following result from \cite[Theorem 1.1]{Har06} (see also \cite[Theorem 4.1]{RT77},  \cite[Theorem 2.13]{Fig17}, \cite[Theorem 1.6.2]{Gut16}). 
\begin{theorem} \label{thm:RT} Let $\nu$ be a positive Borel measure on $D$ such that $\nu(D) < + \infty$ and $h \in  {\rm CV}(D) \cap C^0(\bar D)$.  Then there exists a unique $u  \in {\rm CV}(D) \cap C^0(\bar D)$ such that $\text{M}_\R(v) = \nu $ and $v=h$ in $\partial D$. 
\end{theorem}

The condition $h \in  {\rm CV}(D) \cap C^0(\bar D)$ can be interpreted as the boundary value  being extended as a convex continuous function over $\bar D$. 

\subsection{Relation with the Cegrell class $\mathcal N$}
\begin{proposition}\label{prop: class N continuous}
	Assume $u\in {\rm CV}(D)$ and $u\circ L \in \mathcal N(\Omega)$. Then $u\in {\rm CV}_0(D)$.   
\end{proposition}
\begin{proof}
We fix $w\in \partial D$ and prove that  $\lim_{x\to w} u(x) =0$. Fix a neighborhood $V$ around $w$ and set $U:= V\cap D$. Fix $x_0\in D$ and set 
\[
K:= \frac{1}{2} \bar U + \frac{1}{2} x_0. 
\]
Then $K$ is a compact subset of $D$. 

	Let $(D_j)$ be an increasing sequence of compact subsets of $D$ such that $\cup_j D_j=D$, and set $\Omega_j= L^{-1}(D_j)$. For each $j$ there exists a unique $\varphi_j \in \mathcal F(\Omega)$ such that $(dd^c \varphi_j)^n={\bf 1}_{\Omega_j} (dd^c u_L)^n$. By uniqueness, $\varphi_j=u_j\circ L$ and by Theorem \ref{thm: MA convex} $u_j\in {\rm CV}_0(D)$.  
	There is also a unique $\psi_j \in \mathcal N(\Omega)$ such that 
	\[
	(dd^c \psi_j)^n = {\bf 1}_{\Omega\setminus \Omega_j} (dd^c u_L)^n, \; u_L=u\circ L. 
	\]
	Then $\psi_j=v_j\circ L$ increases to $0$ as $j\to +\infty$. The convergence $v_j \to 0$ is thus uniform on $K$, since they are all continuous on the compact set $K$. By the comparison principle we also have $\psi_j +\varphi_j \leq u_L$, hence $u_j+v_j \leq u$. For $y\in U$, we have, by convexity of $v_j\leq 0$, 
	\[
	v_j(y) \geq v_j(y) + v_j(x_0) \geq  2 v_j\left(\frac{y+x_0}{2}\right)  \geq 2\inf_K v_j. 
	\]
	Thus, for fixed $j$, 
	\[
	\liminf_{y\to w} u(y) \geq \lim_{y\to w} (u_j(y) + v_j(y))  \geq \liminf_{y\to w} v_j(y) \geq 2  \inf_K v_j,
	\]
	using that $u_j \in {\rm CV}_0(D)$. Letting $j\to +\infty$, we arrive at the result. 
\end{proof}

\begin{theorem}\label{thm: MA convex N}
	Assume $\nu$ is a positive Borel measure on $D$ with positive mass $\nu(D)>0$. If there exists $v \in {\rm CV}(D)$ such that $v<0$ and $\int_D (-v) d\nu <+\infty$, then the Monge-Amp\`ere equation ${\rm M}_{\mathbb R}(u) =\nu$ has a unique solution $u\in {\rm CV}_0(D)$.  
\end{theorem}
\begin{proof}
	We write $\psi = v\circ L$ and set $\nu_j = {\bf 1}_{D_j} \nu$, where $D_1\Subset D_2 \cdots \Subset D$ and $\cup_j D_j =D$. By Theorem \ref{thm:RT}, there exists $u_j \in {\rm CV}_0(D)$ such that ${\rm M}_{\mathbb R}(u_j)= \nu_j$. If we write $\varphi_j = u_j \circ L$. then $(dd^c \varphi_j)^n \leq (dd^c \varphi_{j+1})^n$. Indeed, since these functions are toric, it suffices to check the inequality over toric test functions: 
	\[
	\int_{\Omega} \chi\circ L (dd^c \varphi_j)^n = \int_{D} \chi  {\rm M}_{\mathbb R}(u_j) \leq \int_{D} \chi  {\rm M}_{\mathbb R}(u_{j+1}) = \int_{\Omega} \chi\circ L (dd^c \varphi_{j+1})^n,
	\]
	as follows from the formula \eqref{eq:MAR-MAC}. 
	Hence by the comparison principle $\varphi_j$ is decreasing in $j$. Moreover, 
	\[
	\int_{\Omega} (-v\circ L) (dd^c \varphi_j)^n = \int_D (-v) {\rm M}_{\mathbb R}(u_j) \leq \int_D (-v) d\nu<+\infty.  
	\]
	By \cite[Proposition 5.1]{Ceg08}, $\varphi:= \lim_j \varphi_j$ is in $\mathcal E(\Omega)$. If we write $\varphi = u\circ L$, with $u\in {\rm CV}(D)$, then ${\rm M}_{\mathbb R}(u)=\nu$. On the other hand, by \cite[Proposition 5.2]{Ceg08} there exists $\tilde{\varphi} \in \mathcal N(\Omega)$ such that $(dd^c \tilde{\varphi})^n= (dd^c \varphi)^n$. By the comparison principle, $\tilde{\varphi}\leq \varphi_j$ for all $j$, hence $\tilde{\varphi}\leq \varphi$, thus $\varphi\in \mathcal N(\Omega)$, and the comparison principle again gives $\varphi=\tilde{\varphi}$. By Proposition \ref{prop: class N continuous} we thus have $u\in {\rm CV}_0(D)$. 
\end{proof}

\subsection{The eigenvalue problem }
 Let $\nu$ be a positive Borel measure with positive mass on $D$. 
The eigenvalue problem for the real Monge-Amp\`ere operator consists of finding a pair $(\lambda,u)$, where $\lambda >0$ and $u\in \mathrm{CV}_0(D)$ such that
\begin{equation} \label{eq:eigenvalueR}
\text{M}_\R(u) = (-\lambda u)^n \nu, \;  u<0 \; \text{in}\; D.
\end{equation}
This problem has been considered by Lions in \cite{Lio85} for strictly convex domains with smooth boundary and $\nu$ the Lebesgue measure on $\R^n$. He then proved the existence of a unique $\lambda > 0$ and a 
unique (up to a multiple positive constant) smooth convex function $u$ such that $(\lambda,u)$ is a solution to the problem \eqref{eq:eigenvalueR}. For general convex domains and $\nu=dV$ the problem was studied by L\^e in \cite{LeQN18}. 

\smallskip

We will extend these results and prove the following theorem.
 \begin{theorem} \label{thm:eigenvalueR} {\; }
 \begin{enumerate}
 	\item There exists at most one solution to the eigenvalue value problem \eqref{eq:eigenvalueR}. 
 	\item If $\nu =  \text{M}_\R(v)$ for some $v\in {\rm CV}(D) \cap C^0(\bar D)$, then  the eigenvalue problem has a unique solution.  Moreover for any $u_0 \in {\rm CV}_0(D)$ such that $u_0<0$ in $D$,  the iterative sequence defined by
 $$
 \text{M}_{\R} (u_{k+1}) = R(u_k) (-u_k)^n \nu, \, \, \, k \in \N, \; u_{k+1} \in {\rm CV}_0(D),
 $$
 provides a sequence $(u_k)$ in $ {\rm CV}_0(D)$ such that the sequence $(u_k,R(u_k))_k$ converges to $(\lambda_1,u_{\infty})$ where $u_{\infty} \in {\rm CV}_0(D)$ is an eigenfunction of the Monge-Amp\`ere  operator $ \text{M}_\R$.
 \end{enumerate}
  \end{theorem}
  \begin{remark}
  	If $\int_D (-v) d\nu<+\infty$ for some negative convex function $v\not \equiv 0$, then, by Theorem \ref{thm: MA convex N}, there exists $u \in {\rm CV}_0(D)$ such that  $\nu = \text{M}_\R(u)$. Therefore, Theorem \ref{thm:eigenvalueR} applies for any Borel measure with this integrability property. In particular, it applies for any Borel measure with positive finite total mass. 
  \end{remark}
  
  \begin{remark}
  	We can define  $\mathcal E^1(D)$ to be the class of functions $u\in {\rm CV}_0(D)$ such that 
  	\[
  	E_{\mathbb R}(u):=\int_D (-u) {\rm M}_{\mathbb R}(u)<+\infty,
  	\] 
  	and $\mathcal E^1_C(D)$ to be the class of all functions $u$ in $\mathcal E^1(D)$ such that $E_{\mathbb R}(u)\leq C$. Then the condition $\nu = {\rm M}_{\mathbb R}(v)$, for some $v\in {\rm CV}(D)\cap C^0(\bar D)$, can be  replaced by the requirement that $I_{\nu}$ is continuous on each $\mathcal E^1_C(D)$.
  \end{remark}
 \begin{proof} 
 The first statement follows from the uniqueness of solution to the complex Monge-Amp\`ere eigenvalue problem (see Corollary \ref{cor: unique}, Theorem \ref{thm: subsol is sol} and Theorem \ref{thm: uniqueness of the eigenvalue}).
   
To prove the existence part, let $\mu_L := (dd^c v_L)^n$, where $v_L := v \circ L$, and observe that $v_L\in \mathcal E^0(\Omega) \cap C^0(\bar \Omega)$. 
  Then $\mu_L $ is a toric (locally finite) Borel measure on $\Omega$ such that $L_*(\mu_L) = \nu$. By Proposition \ref{prop: I cont},  $I_{\mu_L}$ is continuous on $\mathcal E^1_C$ for any $C > 0$. Therefore, by Theorem \ref{thm: existence} there exists a unique number $\lambda >0$ and  a unique normalized $\varphi \in \mathcal E^1(\Omega)$ such that $(\lambda,\varphi)$ is a  solution to the following problem
  \begin{equation} \label{eq:eigenvalueC}
  (dd^c \varphi)^n = (-\lambda \varphi)^n \mu_L.
  \end{equation} 
  Then  by uniqueness $\varphi$ is a toric plurisubharmonic function in $\Omega$.
  Hence there exists a convex function $u$ in $D$ such that $\varphi = u\circ L$.
  Therefore by \eqref{eq:MAR-MAC} it follows that 
  $$
  \text{M}_\R(u) = (-\lambda u)^n \nu.
  $$
  Moreover since $\int_\Omega (-\varphi)^{n+1} d \mu_L < +\infty$, it follows that 
  $\int_D (-u)^{n+1} d \nu < +\infty$. Hence the Borel measure $\tilde \nu := (-\lambda u)^n \nu$ is of finite mass since $\nu$ has finite mass. 
  By Theorem \ref{thm:RT} there exists a unique $\tilde u \in {\rm CV}_0(D)$ such that 
   \[
   \text{M}_\R(\tilde u) = \tilde \nu := (-\lambda u)^n \nu =  \text{M}_\R(u).
   \]
   Since $\tilde{\varphi} = \tilde{u}\circ L \in \mathcal E^1(\Omega)$ and $(dd^c \varphi)^n= (dd^c \tilde{\varphi})^n$, we infer by uniqueness (see \cite{Ceg98}) that $\tilde{\varphi}=\varphi$, hence $u = \tilde u \in {\rm CV}_0(D)$. This proves the existence part of the theorem.  
   
   The last assertion follows from Theorem \ref{thm: existence iterative}.
  \end{proof}
  
  \begin{example}
  	\label{example: mass}
  	Let $u\in {\rm CV}_0(D)$, $u\neq 0$, and $\alpha\in (0,1)$. Then, for $u_{\alpha}:= -(-u)^{\alpha}$, we have $\int_D {\rm M}_{\mathbb R}(u) =+\infty$. 
  \end{example}
 \begin{proof}
 	It suffices to prove the result for the complex Monge-Amp\`ere operator: $\int_{\Omega} (dd^c \varphi_{\alpha})^n=+\infty$, where $\varphi_{\alpha}:= u_{\alpha}\circ L$. For $t>1$, consider the function 
 	\[
 	v_t:= \max(\varphi_{\alpha}, t\varphi).
 	\]
 	Then $v_t=t\varphi$ near the boundary $\partial \Omega$ because $\varphi(z)\to 0$ as $z\to \partial \Omega$, and $0<\alpha<1$. It follows from \cite[Lemma 4.1]{Ceg98} that 
 	\[
 	\int_{\Omega} (dd^c v_t)^n = \int_{\Omega} (dd^c t \varphi)^n = t^n \int_{\Omega} (dd^c \varphi)^n. 
 	\] 
 	If $\int_{\Omega} (dd^c \varphi_{\alpha})^n<+\infty$, then \cite[Lemma 4.2]{Ceg98} yields 
 	\[
 	\int_{\Omega} (dd^c \varphi_{\alpha})^n \geq \int_{\Omega} (dd^c v_t)^n \to +\infty, 
 	\]
 	a contradiction. 
 	Therefore, $\int_{\Omega} (dd^c \varphi_{\alpha})^n=+\infty$ as desired. 
 \end{proof}

 As in Theorem \ref{thm:Existence-Uniqueness2}, we can also study the following Dirichlet problem with continuous boundary values. Given $h \in C^0(\bar D)\cap {\rm CV}(D)$, we want to find $u \in {\rm CV}(D) \cap C^0(\bar D)$ such that
\begin{equation} \label{eq:realMA}
M_\R(u) = F(\cdot,u) \nu, \; \; \text{and} \; \;  u_{\mid \partial D} = h.
\end{equation}
Setting $\varphi := u\circ L$, then the corresponding Dirichlet problem for the complex Monge-Amp\`ere  operator is as follows:
 \begin{equation} \label{complexMA}
(dd^c \varphi)^n  = F_L(\cdot, \varphi) \nu_L, \; \; \text{and} \; \;  \varphi_{\mid \partial \Omega} = h\circ L,
\end{equation}
where $F_L(z,t) := F(L(z), t)$ is a toric function in $z \in \Omega$ and $h_L = h\circ L$ is a toric function on $\partial \Omega$.

Conversely assume that  \eqref{complexMA} admits a unique solution. Then since $\nu_L$ and $F_L$ are toric, $\varphi$ will be a toric plurisubharmonic function in $\Omega$. Therefore it can be written as $\varphi = u\circ L$, where $u$ is a convex function which will be the unique solution to \eqref{eq:realMA}.

Therefore Theorem \ref{thm:Existence-Uniqueness2}  implies the following theorem.

\begin{theorem}  \label{thm:Existence-Uniqueness3} Assume there exists $ \lambda_0 > 0 $ such that
$$
\frac{\partial F}{\partial t} \geq - \lambda_0 , \, \, \mathrm{in} \, \, D \times \R.
$$ 
Then  

(1) if $\lambda_0 < \lambda_1$,  the maximum principle holds for the operator 
$\text{M}_\R - F(\cdot,u) \, \nu$ on ${\rm CV}(D)$  i.e. if $u, v  \in  {\rm CV}(D) \cap C^0(\bar D)$ is a supersolution and $v \in {\rm CV}(D) \cap C^0(\bar D)$ is a subsolution such that $u \geq v$ in $\partial D$, then $u \geq v$ in $D$. 

Moreover the Dirichlet problem \eqref{eq:realMA} admits a unique solution $u  \in {\rm CV}(D) \cap C^0(\bar D)$. 

(2) If $\lambda_0 = \lambda_1$ and $u, v \in {\rm CV}(D) \cap C^0(\bar D)$ are two solutions to \eqref{eq:realMA}  such that $u(x_0) < v(x_0)$ for some point $x_0 \in D$, then the convex envelope $w$ of the function $\inf \{u-v,0\}$ is an eigenfunction associated  to the eigenvalue $\lambda_1$.

\end{theorem}


\end{document}